\documentclass[12pt,a4paper,reqno]{amsart}

\usepackage[pdftex]{graphics}

\usepackage{arcs}

\theoremstyle{plain}
\newtheorem{theorem}{Theorem}[section]
\newtheorem{corollary}[theorem]{Corollary}

\newtheorem{lemma}[theorem]{Lemma}
\newtheorem{proposition}[theorem]{Proposition}
\theoremstyle{definition}

\theoremstyle{remark}

\newtheorem*{remark}{Remark}

\numberwithin{equation}{section}
\numberwithin{theorem}{section}
\numberwithin{table}{section}
\numberwithin{figure}{section}

\newcommand{\C}{\mathbb C}
\newcommand{\R}{\mathbb R}
\newcommand{\Pabn}{P_n^{(a,b)}}
\newcommand{\UN}{\mathrm{U}(N)}
\newcommand{\abs}[1]{\left \lvert #1 \right \rvert}

\renewcommand{\Re}{\operatorname{Re}}
\renewcommand{\Im}{\operatorname{Im}}

\def\({\left(}
\def\){\right)}

\begin{document}
\title[Roots of the derivative of the zeta function]
{Roots of the derivative of the Riemann zeta function\\
  and of characteristic polynomials}

\author[Due\~nez]{Eduardo Due\~nez}
\address{Department of Mathematics\\
    The University of Texas At San Antonio}
\email{eduenez@math.utsa.edu}

\author[Farmer]{David~W.~Farmer}
\address{
    American Institute of Mathematics }
\email{farmer@aimath.org}

\author[Froehlich]{Sara~Froehlich}
\address{
   Department of Mathematics and Statistics\\
   McGill University}

\author[Hughes]{Chris~Hughes}
\address{Department of Mathematics \\
    University of York}
\email{ch540@york.ac.uk}

\author[Mezzadri]{Francesco~Mezzadri}
\address{
    Department of Mathematics\\
    University of Bristol}
\email{f.mezzadri@bristol.ac.uk}

\author[Phan]{Toan~Phan}
\address{
Department of Economics\\
Northwestern University}

\thanks{
Research supported by the
American Institute of Mathematics and the National Science Foundation}

\begin{abstract}
We investigate the horizontal distribution of zeros of the derivative
of the Riemann zeta function and compare this to the radial distribution
of zeros of the derivative of the characteristic polynomial of a
random unitary matrix.  Both cases show a surprising bimodal distribution
which has yet to be explained.
We show by example that the bimodality is a general phenomenon.
For the unitary matrix case we prove a conjecture
of Mezzadri concerning the leading order behavior, and we show that the
same follows from the random matrix conjectures for the zeros of
the zeta function.
\end{abstract}

\maketitle

\section{Introduction}
The zeros of the derivative $\zeta'(s)$ of the Riemann
zeta-function are intimately connected with the behavior of the
zeros of $\zeta(s)$ itself.  Indeed, a theorem by
Speiser~\cite{Spe34} states that the Riemann Hypothesis (RH) is
equivalent to $\zeta'(s)$ having no zeros to the left of the
critical line.  Thus, understanding of the properties the zeros of
$\zeta'(s)$ can provide important tools and insight into the study
of RH. After Speiser's article this idea was explored by
Berndt~\cite{Ber70} and Spira~\cite{Spi65}, but not much progress
was achieved until the work of Levinson and
Montgomery~\cite{LM74}, who proved a quantitative refinement of
Speiser's theorem. They showed that $\zeta(s)$ and $\zeta'(s)$
have essentially the same number of zeros to the left of the
critical line $\sigma = \Re(s)=\frac12$, and proved that as $T
\rightarrow \infty$, where $T$ is the height on the critical line,
a positive proportion of the zeros of $\zeta'(s)$ lie in the
region
\begin{equation}
  \label{eq:lm_res}
  \sigma < \frac12 + (1 + \epsilon) \frac{\log \log T}{\log T},
  \quad \epsilon > 0.
\end{equation}

Consider the group of unitary matrices $\UN$ with probably
distribution given by Haar measure, which is the unique measure
invariant under the left and right action of $\UN$ on itself.
Such a probability space is often known as the \textit{Circular
Unitary Ensemble} (CUE). Let $\Lambda(z)$ be the characteristic
polynomial of a matrix in the CUE. In recent years evidence has
been accumulated suggesting that, in the limit as $T \to
\infty$, the local statistical properties of $\zeta(s)$ can be
modeled by the characteristic polynomials of matrices in the
CUE where $N\approx \log(T/2\pi)$. The connection between the
Riemann zeta-function and characteristic polynomials is
extensive; examples include the distribution of the zeros of
$\zeta(s)$, its value distribution and its moments. (For a
series of review articles on the subject see~\cite{MS05} and
references therein.)

Assuming that random matrix theory (RMT) provides an accurate
description of $\zeta(s)$, the horizontal distribution of the
zeros of $\zeta'(s)$ in proximity of the critical line should
be the same as the radial distribution of the roots of
$\Lambda'(z)$ close to the unit circle.  This idea was first
developed by Mezzadri~\cite{Mez03}, who determined the
distribution of the zeros of $\Lambda'(z)$ that are very far
from the unit circle and conjectured the leading order term of
the distribution very close to the unit circle.  In this paper
we prove his conjecture. We also perform an analogous
calculation for the Riemann zeta-function and conjecturally
find that the result agrees with the RMT model. In addition, we
do numerical computations in both cases and find a surprising
feature in the distribution of zeros of the derivative, namely
that the probability distribution is bimodal.

\section{The zeros of  $\zeta'(s)$.}
We have mentioned that the main motivation for studying the zeros
of the derivative of the Riemann zeta-function is its connection
with RH. Indeed, Levinson and Montgomery's result is the basis for
Levinson's method~\cite{Lev74}, which Conrey~\cite{Con89} used to
prove that at least 40\% of the zeros of the zeta-function are on
the line $\sigma = \frac12$.

Levinson's method involves estimating a weighted average of the zeros
of $\zeta'(s)$ to the left of $\frac12+a/\log T$ for some fixed~$a>0$.
Thus, zeros of $\zeta'(s)$ in the region $\frac12 \le \sigma <
\frac12+a/\log T$ are an inherent loss in Levinson's method.  It would
be useful to understand the magnitude of this loss. Alternatively, if
we could find a lower bound for the number of zeros of $\zeta'$ in this region we
could improve the estimate for the number of zeros on the critical
line.

\begin{figure}[htp]
\begin{center}
  \scalebox{1.00}[1.00]{\includegraphics{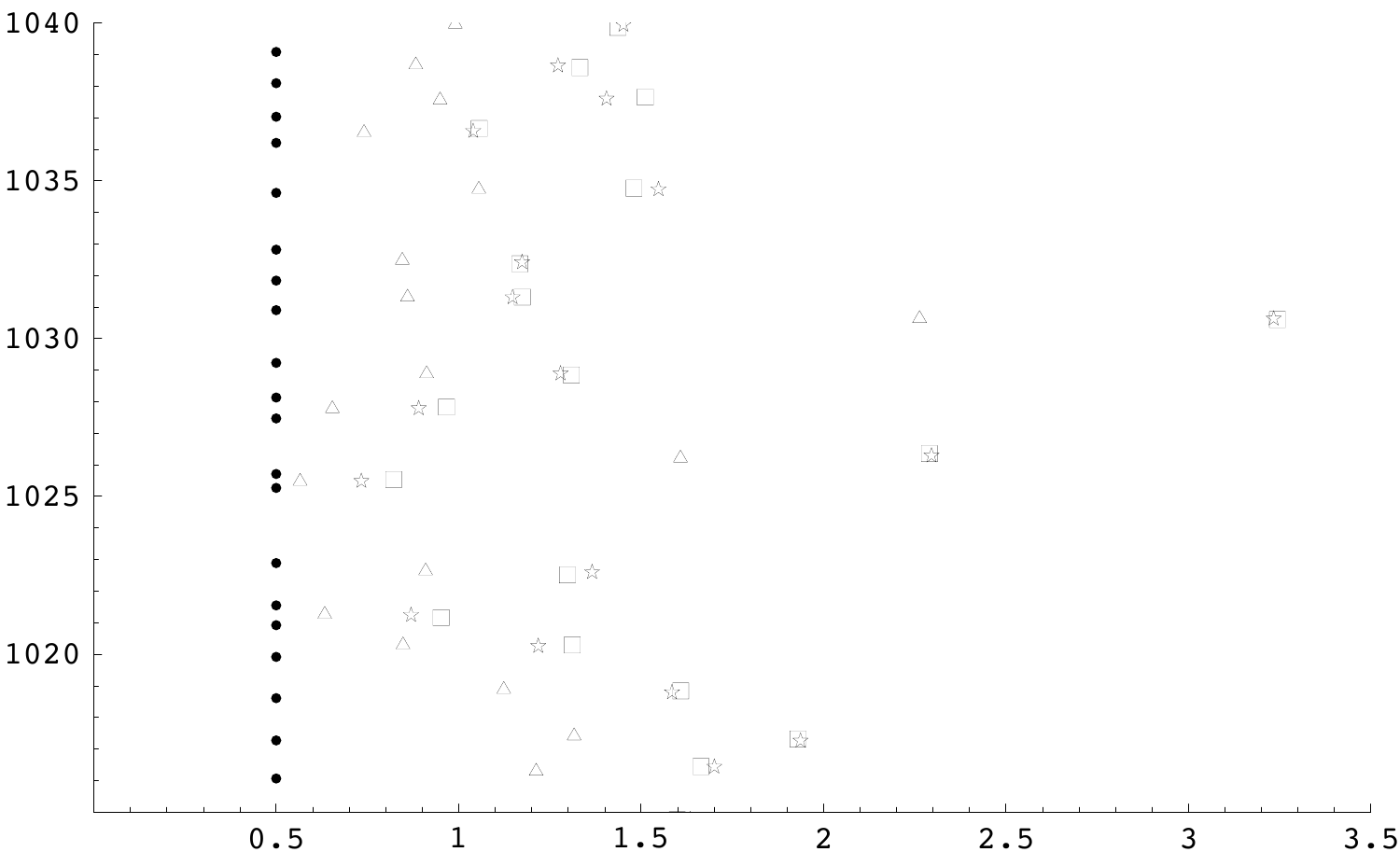}}
  \caption{\sf Zeros of $\zeta(s)$ (dot), $\zeta'(s)$ (triangle),
    $\zeta''(s)$ (square), and $(\zeta'/\zeta)'(s)$ (star), with
    imaginary parts in the range $1015<T<1040$.
\label{fig:zzpzpp}
 }
\end{center}
\end{figure}

Figure~\ref{fig:zzpzpp} gives a representative example of the
location of zeros of $\zeta'(s)$ and the relationship to the zeros
of $\zeta(s)$ and various other derivatives of the zeta-function.
It illustrates that zeros of $\zeta'(s)$ close to the
critical-line correspond to closely spaced zeros of~$\zeta(s)$. We
make this statement precise in Section~\ref{sec:riemann}. Also,
a zero of $\zeta'(s)$ seems to be ``missing'' when zeros of
$\zeta$ are particularly far apart or when there are two
successive large gaps. Indeed, there can't be a zero of
$\zeta'(s)$ between every pair of zeros of $\zeta(s)$ because the
density of zeros of $\zeta(s)$ is $\frac{1}{2\pi}\log(T/2\pi)$
while the density of zeros of $\zeta'(s)$ is
$\frac{1}{2\pi}\log(T/4\pi)$.  So on average there is a
``missing'' zero of $\zeta'(s)$ in each $T$ interval of
width~$2\pi/\log 2 \approx 9.06$.

Conrey and Ghosh~\cite{CG90} and subsequently Guo~\cite{Guo95}
improved Levinson and Montgomery's result~\eqref{eq:lm_res} and showed
that a positive proportion of the zeros of $\zeta'(s)$ are much closer
to the line $\sigma=\frac12$. Indeed, for any fixed $a>0$, the region
\begin{equation}
  \label{eq:cg_res}
  \sigma - \frac12 \ge \frac{a}{\log T}
\end{equation}
contains a positive proportion of the zeros.
Soundararajan~\cite{Sou98} made further progress and introduced the
functions
\begin{subequations}
\label{eq:sound_fns}
\begin{align}
\label{eq:sound_fn-}
  m^-(a)& := \liminf_{T \to \infty} \frac{1}{N_1(T)}
\sum_{\substack{\beta' \le \frac12 + \frac{a}{\log T}\\
                 0 < \gamma' \le T}} 1,\\
\label{eq:sound_fn+}
  m^+(a)& := \limsup_{T \to \infty} \frac{1}{N_1(T)}
\sum_{\substack{\beta' \le \frac12 + \frac{a}{\log T}\\
                 0 < \gamma' \le T}} 1,
\end{align}
\end{subequations}
where $N_1(T)$ is the number of zeros $\beta' + i\gamma'$ of
$\zeta'(s)$ with $0 < \gamma' \le T$.  Soundararajan proved
that $m^-(a)>0$ for $a > 2.6$, and conjectured that
\begin{equation}
  \label{eq:sound_conj}
   m(a) = m^-(a)=m^+(a).
\end{equation}
He also conjectured that $m(a)$ is continuous, that $m(a) > 0$
for all $a >0$, and that $m(a) \to 1$ as $a \to \infty$.
Zhang~\cite{Zha01}, Feng~\cite{Fen05}, Garaev and
Y{\i}ld{\i}r{\i}m~\cite{GY06}, and Ki~\cite{Has06} proved
refinements of Soundararajan's results.  In particular
Feng~\cite{Fen05} showed that $m^-(a)>0$ for all $a>0$
unconditionally of RH but assuming a conjecture on the
frequency of small gaps between consecutive critical zeros of
$\zeta(s)$.

\section{Zeros of derivatives of polynomials
and statement of results}

Suppose $f(z)$ is a polynomial with all zeros on the unit
circle. (Eventually, $f$ will be a random polynomial obtained as the
characteristic polynomial of a random unitary matrix.)  The
Gauss-Lucas theorem assures that all the roots of $f'(z)$ lie on or
inside the unit circle (zeros of $f'$ on the circle occur only if $f$
has multiple zeros). If $f(z)$ has two zeros which are very close
together, then $f'(z)$ will have a zero close by. (This is a
consequence of the continuous dependence of the zeros of $f'$ on those
of~$f$. This dependence is actually piecewise analytic as will be
described below.)  The specific location of the nearby zero of $f'(z)$
will depend primarily on how close those two zeros of $f(z)$ are, and
on the general position of the remaining zeros of~$f(z)$. Thus, to
leading order (with respect to the size of the gap), the distribution
of zeros of $f'(z)$ near $|z|=1$ should largely depend on the
distribution of small gaps between zeros of $f(z)$, that is, on the
tail of the nearest-neighbor spacing of zeros of~$f(z)$.  We will make
this idea precise and treat in detail the case where $f(z)=\Lambda(z)$
is the characteristic polynomial of a CUE matrix (a matrix chosen from
the unitary group~$\UN$, uniformly with respect to Haar measure).

Let $z'$ be a root of $\Lambda'(z)$ and define the random variable
\begin{equation}
  \label{eq:dis_un_c}
    S :=  N(1 - \abs{z'}).
\end{equation}
Denote by $Q(s;N)$ the probability density
function (p.d.f.) of $S$.  Mezzadri~\cite{Mez03} showed that the
limit
\begin{equation}
  \label{eq:limit_pdf}
   Q(s) :=  \lim_{N \to \infty} Q(s;N)
\end{equation}
exists, and proved that
\begin{equation}
  \label{eq:m_r1}
  Q(s;N) \sim \frac{1}{s^2}, \quad N \to \infty, \quad s \to \infty,
\end{equation}
with $s=o\left(N\right)$.  He also conjectured that
\begin{equation}
  \label{eq:m_con}
  Q(s) \sim \frac{4}{3 \pi}s^{1/2}, \quad s \to 0.
\end{equation}

Formula~\eqref{eq:m_r1} can be interpreted as the RMT counterpart of
the Levinson-Montgomery bound~\eqref{eq:lm_res} for the roots of
$\zeta'(s)$.  The RMT model of the Riemann-zeta function is based on
the observation that the local correlations of the non-trivial zeros
of $\zeta(s)$ coincide with those of the eigenvalues of matrices
in the CUE.  In order to make this correspondence quantitative, the
densities of the eigenvalues and of the zeros of $\zeta(s)$ must be
made (asymptotically) equal, i.e.,
\begin{equation}
  \label{eq:dens}
  \frac{N}{2\pi} = \frac{1}{2\pi} \log \frac{T}{2\pi}.
\end{equation}
It follows from~\eqref{eq:m_r1} that the expected value of $S$
does not exist---its p.d.f.\ does not decay sufficiently
rapidly.  On the other hand, using~(\ref{eq:m_r1}), the average
of the values of $S$ not exceeding $N$ is
\begin{equation}
  \label{eq:rmeqlm}
  \sim\int_1^N s\cdot\frac{ds}{s^2} = \log N, \qquad N \to \infty.
\end{equation}
Recalling the relation~(\ref{eq:dis_un_c}) between $S$ and $|z'|$, we
conclude that a positive proportion of the roots of $\Lambda'(z)$ must
lie within a distance from the unit circle bounded from above by
\begin{equation}
  \label{eq:rmeqlm2}
  (1+\epsilon)\frac{\log N}{N}
\end{equation}
from the unit circle. Because of~\eqref{eq:dens},
formula~\eqref{eq:rmeqlm2} corresponds to Levinson and
Montgomery's result~\eqref{eq:lm_res}.

Now consider the roots $e^{it_1}, \ldots, e^{it_N}$ (with $-\pi
< t_i \leq \pi$) of the characteristic polynomial~$\Lambda(z)$
of a random unitary matrix distributed with Haar measure. It is
convenient for our purposes to define
\begin{equation}
  \label{eq:resc_eig}
   x_j = \frac{Nt_j}{2\pi}, \qquad j=1,\ldots,N,
\end{equation}
so that, on average, the distance between two consecutive
$x_j$'s is one. The joint probability density function
(j.p.d.f.) of the eigenvalues is given in terms of
$x_1,\ldots,x_N$ as
\begin{equation}
  \label{eq:weyl_int}
  p_2(x_1,\ldots,x_N) := \frac{1}{N^N N!} \prod_{1 \le j < k \le
    N} \abs{e_N(x_k) - e_N(x_j)}^2,
\end{equation}
where we have used the notation $e_N(x):=\exp(2\pi ix/N)$.
Relabeling the indexes $j=1,\ldots,N$, if necessary, we assume
that
\begin{equation}
  \label{eq:asc_ord}
  x_1 \le \ldots \le x_N < x_1+N.
\end{equation}
% (With this restriction, the j.p.d.f.\ loses the factor of $N!$ on the
% denominator of the right-hand side of~(\ref{eq:weyl_int}).)
We also extend the sequence $\{x_j\}$ to be periodic by setting
$x_{N+1}=x_1$, $x_{N+2}=x_2$, etc.  Fix integers $n$ and $j$
with $0\leq n\leq N-2$. Let us denote by $p_2(n;s)$ the
probability density function of $x_{j+n+1}-x_j$.  Since the
j.p.d.f.~(\ref{eq:weyl_int}) is invariant under translations,
which means
\begin{equation}
  \label{eq:tran_inv}
  p_2(x_1 + \alpha,\ldots,x_N + \alpha)= p_2(x_1,\ldots,x_N)
  \quad \text{for all
  $\alpha \in \mathbb{R}$},
\end{equation}
it follows that $p_2(n;s)$ does not depend on~$j$. If $n=0$,
$p_2(s):=p_2(0;s)$ is known as the \textit{spacing distribution}.  It
has an asymptotic expansion in powers of $s$:
\begin{align}
  \label{eq:sp_as}
    p_2(s) = & \left(\frac{1}{3}-\frac{1}{3 N^2}\right) \pi^2 s^2-
    \left(\frac{2}{45}-\frac{1}{9 N^2}
      +\frac{1}{15 N^4}\right) \pi^4 s^4\cr
    & \quad + \left(\frac{1}{315}-\frac{2}{135 N^2}+
	\frac{1}{45 N^4}-\frac{2}{189 N^6}\right)\pi^6 s^6 + O(s^7).
\end{align}
This expansion follows from the pair correlation function
for $U(N)$ (see \cite{MS05}),
\begin{equation}
\frac{\sin^2(\pi y)}{N^2 \sin^2(\pi y/N)},
\end{equation}
and the fact that the pair correlation function and the nearest neighbor
spacing for $U(N)$ agree to order~6.

To describe our main result, suppose that a root of
$\Lambda(z)$ is degenerate (which means $x_{j+1} = x_j$) so
that $z'= \exp\left(2\pi i x_j/N\right)$ will also be a root of
$\Lambda'(z)$. Simple considerations of continuity show that,
if $x_j$ and $x_{j+1}$ are slightly moved apart (say, while
keeping the remaining roots of $f$ fixed), then $z'$ will also
move, but will still be close to the midpoint of the segment
joining $\exp(2\pi i x_j /N)$ to $\exp(2\pi i x_{j+1}/N)$. In
Proposition~\ref{prop:doma-analyt-root} we make this precise,
showing that $z'$ stays close to that midpoint
provided~$x_{j+1}-x_j<1/\pi$. Henceforth we assume that the
rescaled distance
\begin{equation}
  \label{eq:def_theta}
  \theta: = x_{j+1}-x_j
\end{equation}
is small.

By the translation invariance of the j.p.d.f.\  of the $x_j$'s,
we may assume without loss of generality that
\begin{equation}
  \label{eq:xjxjpu}
  x_{j+1}=\frac{\theta}{2} \quad \text{and} \quad x_j =
  -\frac{\theta}{2}.
\end{equation}
Define
\begin{equation}
  \label{eq:delta_def}
  \delta := N(1 - z').
\end{equation}
In section~\ref{sec:ch_pol_proof} we shall show that
\begin{equation}
  \label{eq:ant_res}
  \delta =  b_1 \pi^2 \theta^2 +  b_2 \pi^4 \theta^4 + O(\theta^6), \quad \text{ as } \theta \to 0,
\end{equation}
where $b_1$ and $b_2$ are explicit functions of the zeros $x_k$ for $k\not= j, j+1$.

By combining the distribution of $\theta$ given in~\eqref{eq:sp_as} with information
we will determine about $b_1$ and $b_2$ in Section~\ref{sec:ch_pol_proof}, we will prove

\begin{theorem}\label{thm:delta}
Let $\Lambda(z)$ be the characteristic polynomial of a random matrix
 in $\UN$ distributed with respect to Haar measure.
The distribution of $\delta=N(1-|z'|)$ arising from closely spaced zeros of
$\Lambda(z)$ is given by
\begin{equation}
\frac{4}{3 \pi} s^{1/2}
- \frac{82}{45\pi}s^{3/2}  + O\left(s^{5/2}\right),
\end{equation}
as $N\to\infty$.
\end{theorem}

Note that Theorem~\ref{thm:delta} refers to the small values of $\delta$ that
arise from closely spaced zeros of the polynomial.
The theorem does \emph{not} account for all small values of~$\delta$.
That distinction is often missed, because
examples such as shown in Figure~\ref{fig:zzpzpp} give the mistaken impression
that zeros of the derivative very close to the unit
circle can only arise from closely spaced zeros of the polynomial.
Farmer and Ki~\cite{FK} give examples of families of polynomials whose
(rescaled) zeros are bounded away from each other,
but for which the density function of $\delta$ vanishes like $C\cdot s$ as $s\to 0$.
They also argue that
any larger density of zeros of $z'$ near the unit circle must arise from
closely spaced zeros of the polynomial.  Therefore we have the following
corollary of Theorem~\ref{thm:delta}, which proves Mezzadri's conjecture
about the distribution of~$|z'|$.

\begin{corollary}
 \label{thm:unitary}
 Let $\Lambda(z)$ be the characteristic polynomial of a random matrix
 in $\UN$ distributed with respect to Haar measure, and let $Q(s;N)$
 be the \textit{p.d.f.} of $S= N(1-|z'|)$ for $z'$ a root of
 $\Lambda'(z)$.  Then
\begin{equation}
Q(s) =  \lim_{N \to \infty} Q(s;N) =  \frac{4}{3 \pi} s^{1/2} +O(s).
\end{equation}
\end{corollary}

Note that it remains an unsolved problem to show that $Q(s)$ is
a proper probability distribution. That is, to show
$\int_0^\infty Q(s)ds=1$. This is the random matrix analogue of
Soundararajan's conjecture~$m(a)\to 1$ as $a\to\infty$.

\section{Comparison with data}
\label{sec:data}

We compare our formulas with numerical data.

We generated Haar-random
matrices in $\UN$ using the simple algorithm
described in \cite{Mezz}.

\begin{figure}[htp]
%\begin{center}
\scalebox{0.7}[0.7]{\includegraphics{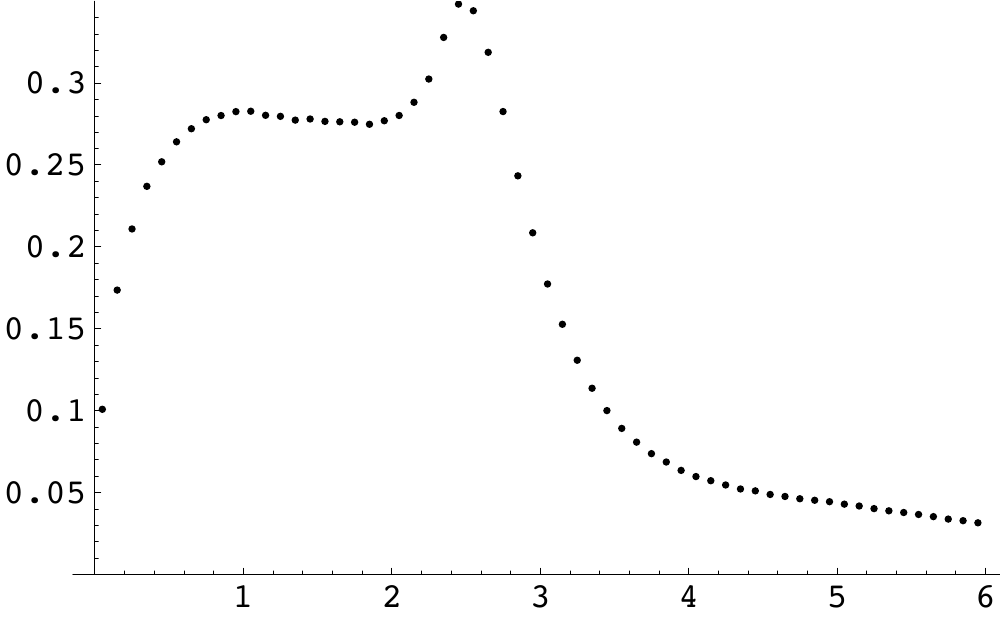}}
\hskip 0.1in
\scalebox{0.7}[0.7]{\includegraphics{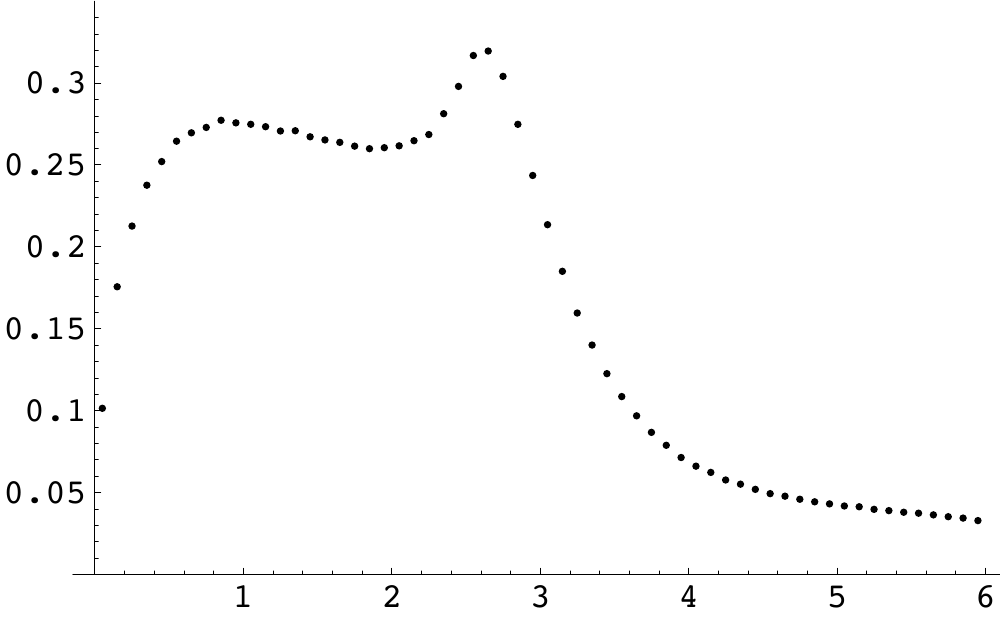}}
\vskip 0.1in
\scalebox{0.7}[0.7]{\includegraphics{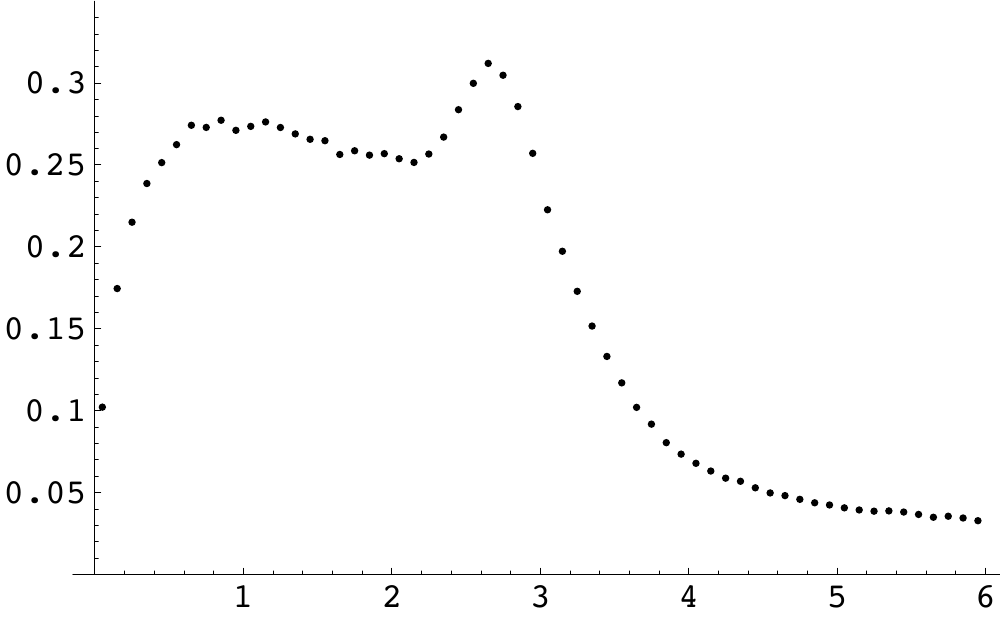}}
\hskip 0.1in
\scalebox{0.7}[0.7]{\includegraphics{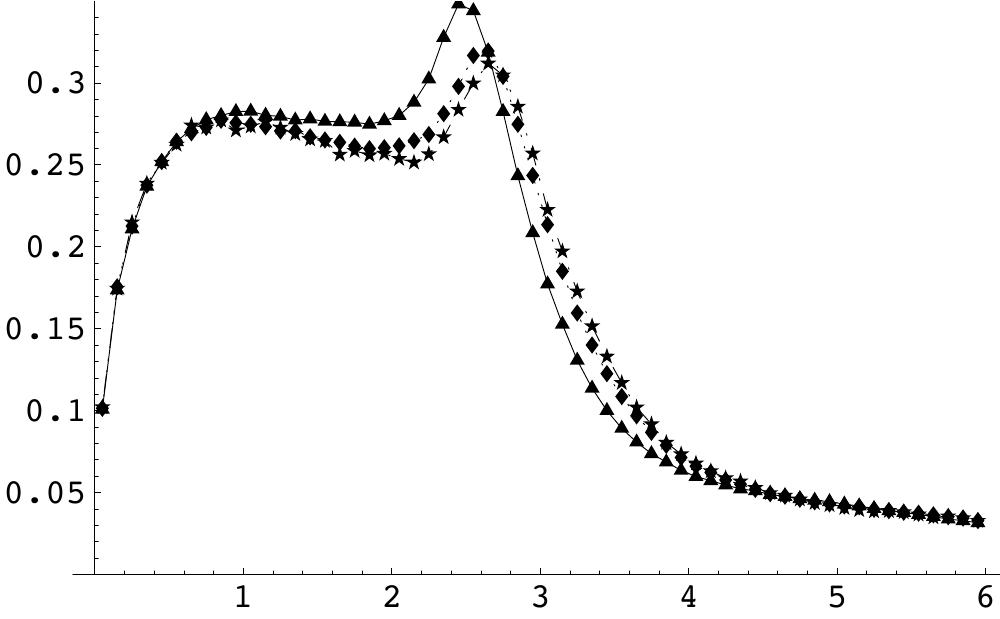}}
\caption{\sf
Distribution of the (normalized) distance from the unit circle of zeros of $\Lambda'(z)$
for $\Lambda$ the characteristic polynomial of a random matrix
from $\UN$. Top row: $N=15$, $40$.  Bottom row: $N=100$
and all three plots together.
} \label{fig:rmtderiv}
%\end{center}
\end{figure}

\begin{figure}[htp]
\begin{center}
\scalebox{0.7}[0.7]{\includegraphics{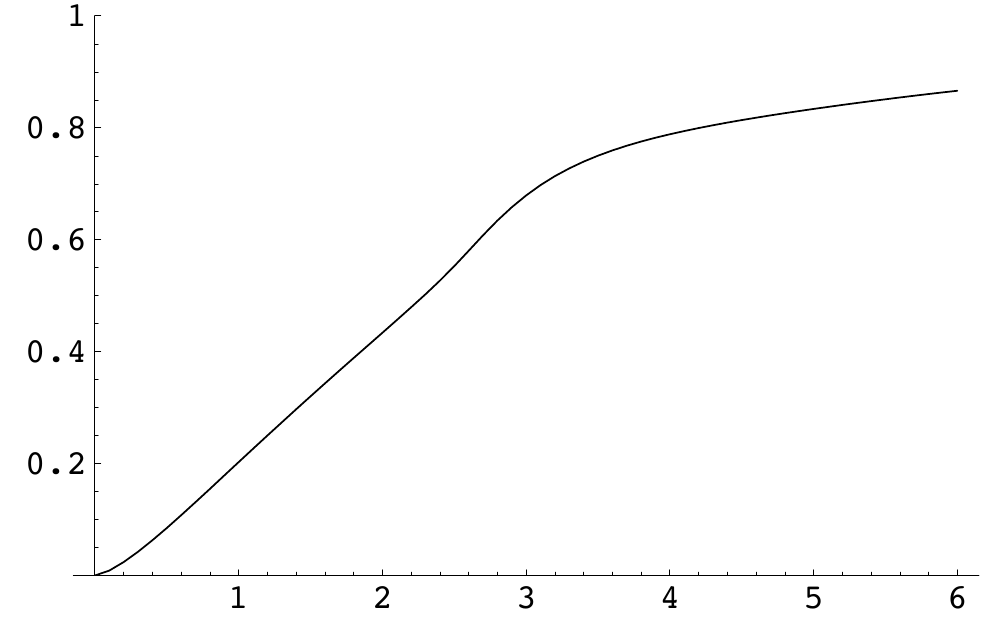}}
\hskip 0.1in
\scalebox{0.7}[0.7]{\includegraphics{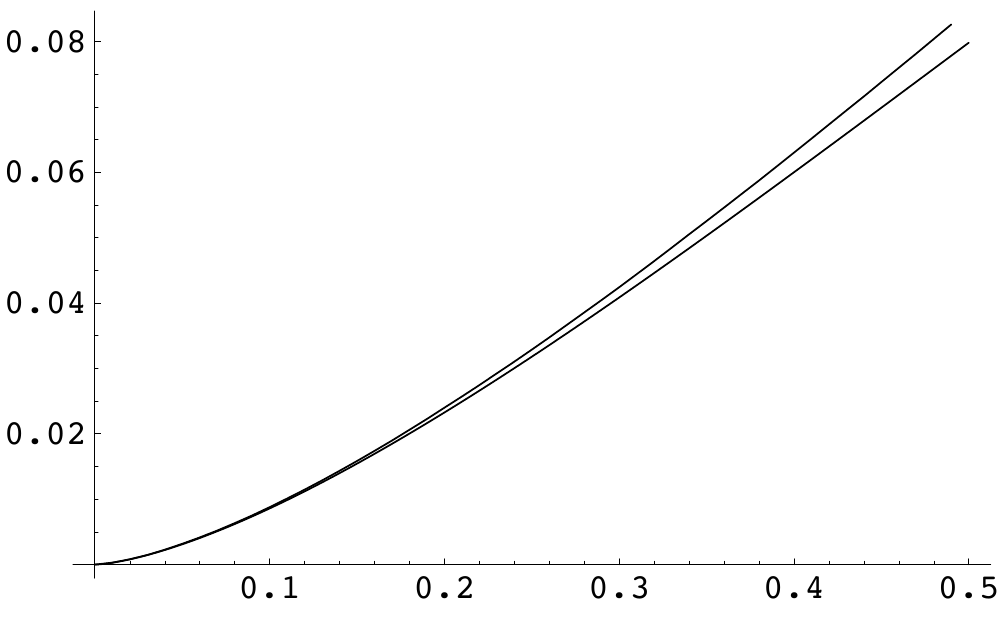}}
\caption{\sf
Experimental cumulative distribution function of $\delta$ for random matrices from  $\mathrm{U}(40)$,
and a comparison of the tail of the distribution with the
main term in Theorem~\ref{thm:unitary}.
} \label{fig:Ip40}
\end{center}
\end{figure}

Figure~\ref{fig:rmtderiv} shows the empirical distribution of the rescaled
zeros of~$\Lambda'$ for various size matrices.
Figure~\ref{fig:Ip40} shows the empirical cumulative distribution
function $Ip(x)=\int_0^x Q(s,40)ds$ for $\mathrm{U}(40)$ and
a comparison with the tail of the empirical cumulative distribution function
with our results, showing good agreement.

\begin{figure}[htp]
%\begin{center}
\scalebox{0.7}[0.7]{\includegraphics{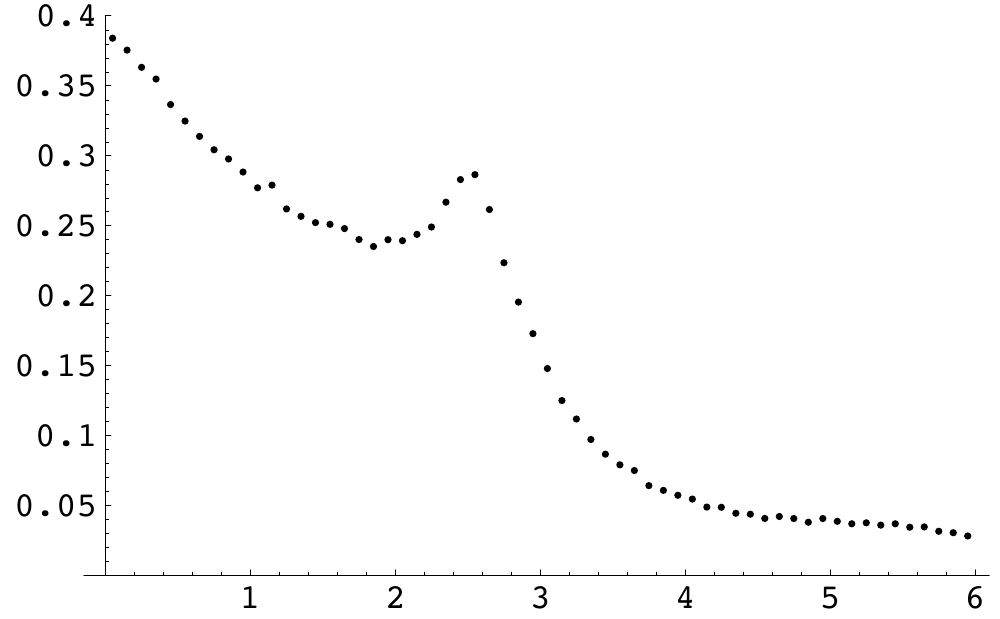}}
\hskip 0.1in
\scalebox{0.7}[0.7]{\includegraphics{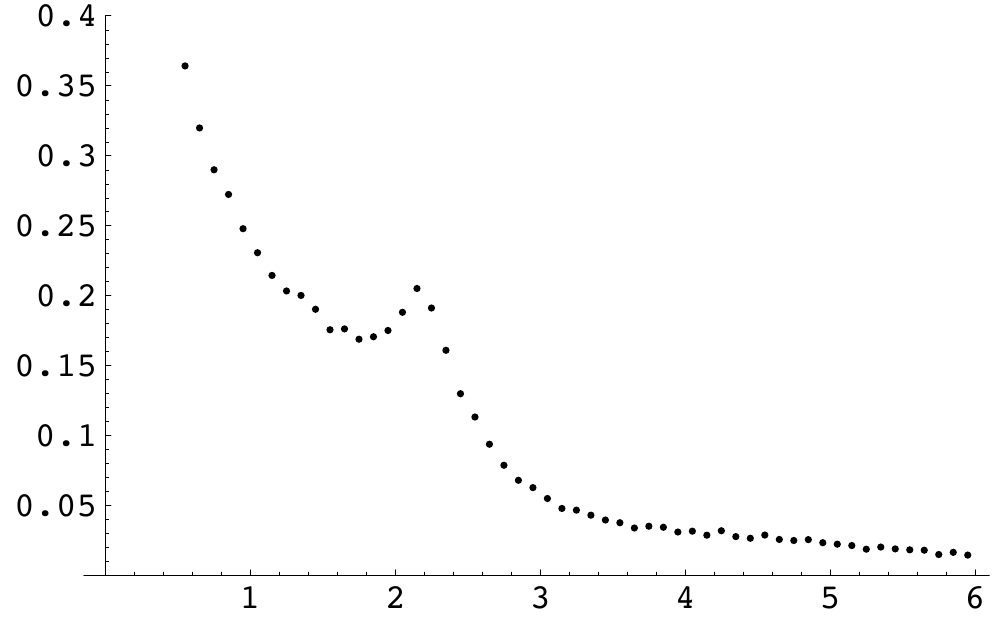}}
\caption{\sf
Distribution of the (normalized) absolute value of zeros of $\Lambda'(z)$
for $\Lambda$ the characteristic polynomial of a random matrix
from $\mathrm{COE}(40)$ (left) and
$\Lambda$ a degree-$40$ polynomial with independent uniformly distributed
(Poisson)
roots on the unit circle (right).  For the Poisson case, the plot is truncated
to suppress the large contribution from the small values.
} \label{fig:coepoiss40}
%\end{center}
\end{figure}

We would like to know the underlying cause of the curious ``second
bump'' in the distribution of zeros of derivatives.  This seems to be
a completely general phenomenon.  In Figure~\ref{fig:coepoiss40} we
show the analogous distributions for characteristic polynomials of
matrices from $\mathrm{COE}(40)$ and for degree-$40$ polynomials whose
roots are independently and uniformly distributed on the unit circle.
Both cases show the ``second bump'', although not quite at the same
location.  In Figure~\ref{fig:rezetaprime} we find a similar shape for
the distribution of zeros of~$\zeta'$.

\section{Proof of theorem~\ref{thm:delta}}
\label{sec:ch_pol_proof}

Suppose $f(z)$ is a degree-$N$ polynomial having all zeros on the unit
circle, for which two zeros $z_1$, $z_2$ are very close together.
Then the derivative $f'(z)$ will have a zero close to the midpoint
$(z_1+z_2)/2$.  This follows because, if $z_1=z_2$ (that is, if $f(z)$
has a multiple root at $z_1$), then $z_1$ is also a root of the
derivative $f'(z)$, and the roots of $f'$ are continuous functions of
the roots of~$f$.  By a rotation we can assume that
\begin{equation}\label{eq:2}
f(z)=F(z)(z-e^{-i\Theta/2})(z-e^{i \Theta/2}),
\end{equation}
where $F(z)$ does not have any zeros $e^{i t}$ with $-\Theta/2\le t\le
\Theta/2$.  The root of $f'$ near~$1$ is the root near~$1$ of
\begin{equation}\label{eqn:fprimeoverf}
\frac{f'}{f}(z)=\frac{F'}{F}(z)
    + \frac{1}{z-e^{-i\Theta/2}} + \frac{1}{z-e^{i \Theta/2}} ,
\end{equation}
and we denote that root by $z'=1-\Delta$.

We are concerned with the case when $f$ is the characteristic
polynomial of a random CUE matrix, and we want to understand the
distribution of the zeros of~$f'$.  Since the CUE measure (i.~e.,
normalized Haar measure on $\UN$) is invariant under
rotation, the distribution depends only on the absolute value of the
roots of~$f'$.  Those roots will accumulate near the unit circle as
$N$ grows, so we must rescale them suitably in order to get a
meaningful result in the limit $N\to\infty$.  We let
\begin{equation}\label{eq:1}
  \begin{aligned}
\Theta &= 2 \pi\theta/N \\
\Delta &= \delta/N .
\end{aligned}
\end{equation}
Note that this rescaling gives $\langle \theta \rangle =1$.
The rescaling of~$\delta$ is more subtle, and it will
be found in equation~\eqref{eqn:deltathetarelationN} that this is the
correct rescaling. Note that, while $\theta$ is a real number,
$\delta$ is usually complex (but typically with small imaginary part,
as we shall see).

We will determine the leading order behavior
of the roots which are close to the unit circle, so in the above
notation we are interested in
\begin{align}\label{eqn:normalized}
N(1-|z'|) =\mathstrut &  N( 1- |1-\Delta| ) \cr
    = \mathstrut & \Re(\delta) - \frac{\Im(\delta)^2}{2N} + O\left(\frac{\delta^3}{N^2}\right) .
\end{align}

\subsection{Expansion for the roots}

Exploiting the symmetry of the \textit{j.p.d.f.}~(\ref{eq:weyl_int})
under arbitrary relabellings of the variables $x_j$, as well as its
translation invariance (cf., Section~\ref{thm:unitary}), we will
assume (without loss of generality) that $F(z)$ from
Equation~(\ref{eq:2}) is given in the form (recall that
$e_N(x):=\exp(2\pi ix/N)$)
\begin{equation}
  F(z) = \prod_{n=1}^{N-2} \left(z-e_N(x_n)\right),\qquad  x_n \in
  \left(-\frac N2,-\frac\theta2\right]\cup\left[\frac\theta2,\frac N2\right).
\label{eq:3}
\end{equation}
(We are not excluding the possibility that $F$ has
roots at $e_N(\pm\theta/2)$.)  Therefore,
\[
\frac{F'}{F}(z) = \sum_{n=1}^{N-2} \frac{1}{z-e_N(x_n)}.
\]
Now let $z=z'$ be a root of $f'$ (say, that root which is closest to
$z=1$).  The Implicit Function Theorem shows that $z'$ is an analytic
function of $\theta$ (at least for $\theta$ sufficiently
small). Indeed, fixing $F$ and regarding $f'$ as a function
$f'(\theta;z)$ of both $\theta$ and~$z$, we have $f'(0;1)=0$; by the
assumption that $F$ has no root at $z=1$ all that remains to observe
is that
\begin{equation}
  \label{eq:not_inline}
  0 \neq 2F(1) = \frac{\partial}{\partial z}f'(0;1)=f''(0;1).
\end{equation}
if $N$ is sufficiently large, by
Proposition~\ref{prop:doma-analyt-root} $z'$ is defined uniquely and
analytically as a function of $\theta$ in the domain
\begin{equation}
  |\theta| < \min\{ \tfrac{1}{\pi}, |x_1|, \ldots, |x_{N-2}|\}.
\end{equation}

We write $z'=1-\delta/N$ and wish to expand
$(1-\delta/N-e_N(x_n))^{-1}$ as a Taylor series in~$\delta$.  This
will be justified when $|\delta| < N |1-e_N(x_n)|$. Hence, we have
\begin{equation}\label{eqn:FprimeoverF}
\frac{F'}{F}\(1-\frac{\delta}{N}\) = N \sum_{j=0}^\infty A_j
\delta^j
\end{equation}
for $\delta$ sufficiently small, where
\begin{align}\label{eqn:Ajdef}
A_j &= (-1)^j\frac{1}{j! N^{j+1}} \left(\frac{F'}{F}\right)^{(j)}(1) \cr
&=\frac{1}{N^{j+1}} \sum_{n=1}^{N-2} \frac{1}{(1-e_N(x_n))^{j+1}}.
\end{align}
We will see later that the prefactor of $N$ is the right choice to
make the coefficients $A_j$ approximately bounded.  Note that, for
$|\theta| < 1$ (equivalently, for $|\Theta|<2\pi/N$), any $\delta$
such that $|\delta|<|\theta|=N|\Theta|/(2\pi)$ makes the expansion in
\eqref{eqn:FprimeoverF} valid (independently of the exact location of
the zeros of~$F$), since
\begin{equation}
  |\delta| < \frac{N \Theta}{2\pi} \leq
  N \left| 1-e_N(\theta/2) \right| \leq
  N \left| 1-e_N(x_n) \right|,\qquad 1\leq n\leq N-2.
\end{equation}

We also wish to expand the other terms in \eqref{eqn:fprimeoverf}
as a series in $\delta$. We have
\begin{multline}
\frac{1}{z'-e^{-i\Theta/2}} + \frac{1}{z'-e^{i \Theta/2}}
=\frac{2-\frac{2\delta}{N} - 2\cos\left(\frac{\pi\theta}{N}\right)}
{2-\frac{2\delta}{N} +\frac{\delta^2}{N^2}
-2\left(1-\frac{\delta}{N}\right) \cos\frac{\pi\theta}{N}} \cr
=\frac{-2\delta+\pi^2 \theta^2 N^{-1} - \frac{1}{12} \theta^4 \pi^4 N^{-3} + O(\theta^6 N^{-5})}
{\delta^2 + \pi^2\theta^2 - \pi^2 \delta \theta^2 N^{-1} - \frac{1}{12}\pi^4\theta^4 N^{-2} + \frac{1}{12}\pi^4\delta \theta^4 N^{-3} + O(\theta^6 N^{-4})}.
\end{multline}
Combining this with equations \eqref{eqn:fprimeoverf} and
\eqref{eqn:FprimeoverF}, putting all terms over a common denominator,
and using the fact that $f'(z')=0$, we have
\begin{multline}\label{eqn:deltathetarelationN}
0= \sum_{j=0}^\infty A_j \delta^j
\left(\delta^2+ \pi^2 \theta^2 - \pi^2 \delta \theta^2 N^{-1}  - \tfrac{1}{12}\pi^4\theta^4 N^{-2} + \tfrac{1}{12}\pi^4\delta \theta^4 N^{-3} + O(\theta^6 N^{-4})\right) \cr
-2\delta + \pi^2 \theta^2 N^{-1} - \tfrac{1}{12} \pi^4 \theta^4 N^{-3} + O(\theta^6 N^{-5}).
\end{multline}
Note that a global factor of $N$ canceled to give the above
equation, suggesting that we have chosen the correct scaling for~$\delta$.

Equation~\eqref{eqn:deltathetarelationN} is simply a more explicit and
manageable form of the equation $f'(\theta;z')=0$ defining $z'$
implicitly as a function of~$\theta$. Noting that $f'(0;1)=0$,
together with the functional equation $f'(\theta;z)=f'(-\theta;z)$, it
follows that $\delta=\delta(\theta)$ has an expansion in powers
of~$\theta^2$, with no constant term, of the form
\begin{equation}\label{eqn:deltaapprox1}
\delta = b_1 \pi^2\theta^2 + b_2 \pi^4 \theta^4 + O(\theta^6).
\end{equation}

From~\eqref{eqn:deltathetarelationN}, we obtain
\begin{multline}\label{eqn:subdelta}
  0= \left( A_0 - 2 b_1 + \frac{1}{N} \right)\pi^2 \theta^2 \cr
  +\left(A_1 b_1 + A_0 b_1^2 - 2 b_2 - \frac{A_0 b_1}{N} -
    \frac{A_0}{12N^2} - \frac{1}{12N^3} \right)\pi^4 \theta^4
  +O(\theta^6).
\end{multline}
Setting each term in~\eqref{eqn:subdelta} equal to 0 and solving for
$b_1$ and $b_2$ we have the following:

\begin{proposition}\label{prop:delta} In the notation above, if $0\le\theta<1/\pi$
and $N$ is sufficiently large then
$\delta= b_1 \pi^2\theta^2 + b_2 \pi^4 \theta^4 + O(\theta^6)$ where
  \begin{align}
\label{eqn:bAN}
  b_1 &= \frac{A_0}{2} + \frac{1}{2N} \\
  b_2 &=\frac{1}{8}\left(A_0^3 +2 A_0 A_1\right) + \frac{A_1}{4N} -
  \frac{A_0}{6N^2} - \frac{1}{24 N^3} ,
\end{align}
with  $A_j$ given in~\eqref{eqn:Ajdef}.
\end{proposition}

Note that in this analysis we have treated $F$ (and hence $A_j$) as
being fixed, in the sense that we assume its zeros do not vary with
$\theta$. In the next section we will show that
when $f(z)$ is the characteristic polynomial of a
matrix drawn from the CUE, this can be justified up to
$O(\theta^7)$.

Using \eqref{eqn:deltaapprox1} and \eqref{eqn:bAN}, we can determine the
distribution of $\delta$ from the distributions of $\theta$
and the~$A_j$.  For small $\theta$, this comes from the tail of the
nearest-neighbor spacing.

\subsection{Nearest-neighbor spacing}
\label{nnsp}

Proposition~\ref{prop:delta}
provides a formula for $\delta=N(1-z')$, but what we really want is
the distribution of $\delta^*:=N(1-|z'|)$.
So by \eqref{eqn:normalized} and \eqref{eqn:deltaapprox1}, and writing
$B_j=\Re(b_j)$, we have
\begin{align}\label{eqn:delta*}
\delta^* =\mathstrut& B_1 \pi^2\theta^2 + B_2 \pi^4 \theta^4
+ O\left(\theta^6 + \frac{|\delta|^2}{N} + \frac{|\delta|^3}{N^2} \) \\
=\mathstrut& B_1 \pi^2\theta^2 + B_2 \pi^4 \theta^4
+ O\left(\theta^6 + \frac{\theta^2}{N}\).
\end{align}
The second line is a corollary of  Proposition~\ref{prop:doma-analyt-root},
because $\delta\ll\theta$ if $\theta<1/\pi$.

Suppose for the moment that
$B_1$ and $B_2$ were
constants (instead of being random).
We would have
$\delta^*=g(\theta)=B_1\pi^2\theta^2 + B_2\pi^4\theta^4 +O_N(\theta^6)$
where $\theta$ is random
with p.d.f.~\eqref{eq:sp_as}
given by the nearest neighbor spacing of eigenvalues of unitary matrices.
Then the distribution function of
$ \delta^*$ would be given by
\begin{align}\label{eqn:deldist}
\frac{p_2(g^{-1}(s))}{g'(g^{-1}(s))}
=\mathstrut&
\frac{B_1^{-3/2}}{6 \pi} \left(1-\frac{1}{N^2}\right)  s^{1/2} \cr
&-
\frac{1}{\pi}
\left(\frac{ B_1^{-5/2}}{45} \left(1-\frac{5}{2N^2} +\frac{3}{2N^4}\right)
 + \frac{5 B_1^{-7/2} B_2}{12}\left(1-\frac{1}{N^2}\right)
 \right) s^{ 3/2} \cr
&+ O_N(s^{5/2}).
\end{align}

It turns out that $B_1$ actually is a constant: in
Section~\ref{sec:Aj} we show that $B_1=\frac{1}{4}$.  It is fortunate
that $B_1$ is a constant, otherwise it could be difficult to determine
the expected value of quantities like~$B_1^{-3/2}$.

The contribution of $B_2$ takes a bit more work.  If $B_2$ was
independent of $\theta$ then we could just average over the possible
contributions of $B_2$ to~\eqref{eqn:deldist}.  That is, in
\eqref{eqn:deldist} replace $B_2$ by its expected value.  This can be
computed from the expected values of various combinations of $A_0$ and
$A_1$.  But $B_2$ is not independent of $\theta$.  However, it is
independent of $\theta$ to leading order.  Our specific concern is the
distribution of the other roots when $\theta$ is very small.  This
approximates the polynomial having a double zero.  Since the
next-nearest-neighbor spacing of $\UN$ eigenvalues vanishes to
order~$7$, the dependence of $B_2$ on $\theta$ is only to order
$O(\theta^7)$.  Thus, for the terms we are computing for the
distribution of $\delta^*$ we can treat $B_2$ as independent
of~$\theta$.  The expected value of $B_2$ is calculated in
Section~\ref{sec:Aj}.

%Note that the nearest-neighbor spacing \eqref{eq:sp_as} is just the first three terms in the
%expansion of the pair correlation function for~$\UN$,
%\begin{equation}
%1 - \left(\frac{\sin(\pi s)}{N \sin( \pi s/N)}\right)^2 .
%\end{equation}
%The next-nearest-neighbor spacing of eigenvalues in $\UN$
%vanishes to order~7, so the nearest neighbor spacing and
%the pair correlation agree to order~6.

\subsection{Expected value of $A_j$}\label{sec:Aj}

We have
\begin{align}\label{eqn:Ajsum}
A_j :=\mathstrut &\frac{1}{N^{j+1}j!} \left(\frac{F'}{F}\right)^{(j)}(1) \cr
=\mathstrut &(-1)^j N^{-j-1} \sum_{n=1}^{N-2} \frac{1}{\left(1-e^{i t_n} \right)^{j+1} },
\end{align}
where $t_1,t_2,\ldots$ are the arguments of the zeros of $F(z)$.
Since
\begin{equation}
\frac{1}{1-e^{i t}} =   \frac12 + \frac{i}{2} \cot\(\frac{t}{2}\)
\end{equation}
we see that
\begin{equation}
\Re(A_0)= \frac{N-2}{2 N},
\end{equation}
so
\begin{equation}
B_1=\Re(b_1)=\frac14 ,
\end{equation}
as claimed.
Note also that
\begin{equation}
\langle A_0\rangle= \frac{1}{2}-\frac{1}{N},
\end{equation}
because the imaginary part of the summand is odd.

For $A_j$ with $j\ge 1$, we require a random matrix calculation.
The sum in \eqref{eqn:Ajsum} is dominated by the terms
where $e^{i t_n}$ is close to~$1$, so one possibility is to determine
the level densities of the~$t_n$.   We will find the expected value
of $A_j$ by appealing to prior results on averages of ratios
of characteristic polynomials~\cite{CFZ}.

We assume that $f(z)$ is the characteristic polynomial of a matrix
chosen uniformly with respect to Haar measure from the unitary group~$\UN$.
We restrict to those matrices which
have two eigenvalues very close to~$1$, and we wish to determine
the joint distribution of the remaining eigenvalues.  This is very
similar to the calculations of Due\~nez~\cite{D} and Snaith~\cite{Sna05}
for the orthogonal group~$\mathrm{SO}(N)$.

First we restrict the measure on the entire ensemble to determine the
measure on the remaining eigenvalues.  Haar measure on $\UN$
is given by
\begin{equation}
d\mu = C \prod_{1\le n<m\le N}
    \left|e^{i t_n} - e^{i t_m}\right|^2 dt_1 \cdots dt_N .
\end{equation}
Here and following, $C$ is a normalization constant which may
vary from line to line, chosen so that the measure has total mass~$1$.
Restricting to those matrices which have \emph{one} eigenvalue
at~$1$ is equivalent to rotating (changing variables) to move
an eigenvalue to~$1$.  So we can also write the measure as
\begin{equation}\label{eqn:measure1}
d\mu_{1} = C
\prod_{1\le n<m\le N-1}
    \left|e^{i t_n} - e^{i t_m}\right|^2
\prod_{1\le n\le N-1}
    \left|e^{i t_n} -1\right|^2
dt_1 \cdots dt_{N-1} .
\end{equation}
The set of matrices which have a repeated eigenvalue at~$1$ has measure zero,
so there is no canonical way to restrict the measure.  However, we are
interested in the limiting case of two eigenvalues which are very
close together, so we determine the measure by restricting the
measure~\eqref{eqn:measure1} to have $|t_{N-1}|\le t$, and then
let $t\to 0$.  The resulting measure is
\begin{equation}\label{eqn:measure2}
d\mu_{2} = C
\prod_{1\le n<m\le N-2}
        \left|e^{i t_n} - e^{i t_m}\right|^2
\prod_{1\le n\le N-2}
        \left|e^{i t_n} -1\right|^{4}
dt_1 \cdots dt_{N-2} .
\end{equation}

Let $U_2(N-2)$ denote the ensemble of unitary matrices with
joint eigenvalue
measure~$\mu_2$.  Then if $g=g(e^{i t_1},\ldots,e^{i t_{n-2}})$ we have
\begin{equation}
\langle g\rangle_{U_2(N-2)} =
C_N
\langle g |\Lambda(1)|^4 \rangle_{\mathrm{U}(N-2)}
\end{equation}
where the right side is a Haar measure average, $\Lambda$
is the characteristic polynomial, and
\begin{equation}
C_N=\langle |\Lambda(1)|^4 \rangle_{\mathrm{U}(N-2)}^{-1} .
\end{equation}
In other words, an expectation involving repeated eigenvalues on $\UN$
is equivalent to an expectation on $\mathrm{U}(N-2)$ with an extra factor of the
$4^{th}$ power of the characteristic polynomial.
This is the key observation for
computing the expected values of the $A_j$ because it reduces it
to the evaluation of known quantities.

Specifically,
let
\begin{equation}
G(\alpha_1,\alpha_2,\alpha_3,\alpha_4; \beta_1,\beta_2; \gamma_1,\gamma_2)
=
\frac{
\Lambda(e^{-\alpha_1})
\Lambda(e^{-\alpha_2})
\Lambda(e^{-\alpha_3})
\Lambda(e^{-\alpha_4})
\overline{\Lambda}(e^{-\beta_1})
\overline{\Lambda}(e^{-\beta_2})
}
{\Lambda(e^{-\gamma_1})
\Lambda(e^{-\gamma_2})
} .
\end{equation}
Theorem 4.1 of \cite{CFZ} provides an explicit formula for the expected
value of
$G$ for $\Lambda$ the characteristic polynomial of
Haar distributed matrices on~$\mathrm{U}(N-2)$.
The formula is complicated so we do not reproduce it here.
This is sufficient to determine
the expected values of all the quantities in~\eqref{eqn:bAN}.
The calculation requires the assistance of a computer algebra
package.  We now present the answers, which we determined with
the help of Mathematica.

The normalization constant for the measure $\mu_2$ is (the reciprocal of)
\begin{align}
\langle G(0,0,0,0;0,0;0,0)\rangle_{\mathrm{U}(N-2)}
&= \langle |\Lambda(1)|^4 \rangle_{\mathrm{U}(N-2)}\cr
&=
\frac{N^4}{12}-\frac{N^2}{12} \cr
&= C_N^{-1},
\end{align}
say.  As $\theta\to 0$ we have
\begin{align}
\langle A_0^3\rangle =\mathstrut&
\left\langle \left(\frac{\Lambda'}{\Lambda}\right)^3 (1)\right\rangle_{U_2(N-2)} \cr
=\mathstrut& - C_N
   \frac{\partial^3}{\partial_{\alpha_1,\alpha_2,\gamma_1}}
\bigg|_{(\alpha_1,\alpha_2,\gamma_1) =(0,0,0)}
\left\langle
G(\alpha_1,\alpha_2,0,0; 0,0; \gamma_1,0) \right\rangle_{\mathrm{U}(N-2)}
\cr
=\mathstrut&
\frac{1}{10} N^3-\frac{7}{10} N^2+\frac{8}{5 }N-\frac{6}{5}\cr
&\cr
\langle A_1 \rangle =\mathstrut&
\left\langle \left(\frac{\Lambda'}{\Lambda}\right)' (1)\right\rangle_{U_2(N-2)}\cr
=\mathstrut& C_N
\left(1+\frac{d}{d\alpha_1}\right)\bigg|_{\alpha_1=0}
\
\frac{\partial}{\partial \alpha_1} \bigg|_{\gamma_1=\alpha_1}
\left\langle
G(\alpha_1,0,0,0; 0,0; \gamma_1,0)
\right\rangle_{\mathrm{U}(N-2)}\cr
=\mathstrut&
\frac{1}{15}N^2-\frac{1}{2} N+\frac{11}{15}  \cr
&\cr
\langle A_0 A_1 \rangle =\mathstrut&
\left\langle \left(\frac{\Lambda'}{\Lambda}\right)(1) \left(\frac{\Lambda'}{\Lambda}\right)'(1) \right\rangle_{U_2(N-2)} \cr
=\mathstrut& -C_N
\left(1+\frac{d}{d\alpha_1}\right)\bigg|_{\alpha_1=0}
\
\frac{\partial}{\partial \alpha_1} \bigg|_{\gamma_1=\alpha_1}
\
\frac{\partial}{\partial \alpha_2} \bigg|_{(\alpha_2,\gamma_2)=(0,0)} \cr
&\phantom{XXXXX}
\left\langle
G(\alpha_1,\alpha_2,0,0; 0,0; \gamma_1,\gamma_2)
\right\rangle_{\mathrm{U}(N-2)}\cr
=\mathstrut&
\frac{1}{30} N^3-\frac{3}{10}N^2+\frac{13}{15}N-\frac{4}{5}.
\end{align}
Thus,
\begin{equation}\label{eqn:b2}
\langle B_2 \rangle=\Re \langle b_2 \rangle =\frac{1}{48}-\frac{7}{48}N^{-1}+ O(N^{-2}).
\end{equation}
Inserting this into
\eqref{eqn:deldist} gives the expansion claimed in Theorem~\ref{thm:delta}.

\section{The Riemann zeta-function}\label{sec:riemann}

We do analogous calculations for derivatives of the Riemann zeta-function
and compare our results with data.

We start with
\begin{equation}
\frac{\zeta'(s)}{\zeta(s)} = b-\frac{1}{s-1}-\frac{1}{2}\frac{\Gamma'(\frac{1}{2}s+1)}{\Gamma(\frac{1}{2}s+1)}+\sum_\rho\left(\frac{1}{s-\rho}+\frac{1}{\rho}\right)
\end{equation}
where $b= \log2\pi-1-\frac{1}{2}\gamma$.
As in the polynomial case, we assume there are two very closely spaced zeros
of the $\zeta$-function and look for the nearby zero of $\zeta'$.
Suppose the closely spaced zeros are
\begin{equation}
\rho_\pm = \frac12 + i (t\pm \tfrac12 \Theta)
\end{equation}
with
\begin{equation}
s'= \frac12 + X + i t
\end{equation}
a zero of~$\zeta'$.

Using
\begin{equation}
\frac{\Gamma'}{\Gamma}(s) = \log(s) + O(1/s)
\end{equation}
we have
\begin{equation}
0= b^* -\frac12 \log t + \frac{1}{s'-\rho_-}+\frac{1}{s'-\rho_+}
  + \sum_{\rho\not= \rho_\pm} \left(\frac{1}{s'-\rho}+\frac{1}{\rho}\right),
\end{equation}
where
\begin{equation}
b^* = b + \frac12\log(2) - i\frac{\pi}{4} + O(1/t).
\end{equation}

Note that
\begin{equation}
\frac{1}{s'-\rho_-}+\frac{1}{s'-\rho_+}= \frac{8 X}{4 X^2 +\Theta^2}
\end{equation}
and
\begin{equation}
\frac{1}{s'-\rho}+\frac{1}{\rho}
= \frac{X-i(t-\gamma)}{X^2+(t-\gamma)^2} + \frac{\frac12 -i \gamma}{\frac14 + \gamma^2} ,
\end{equation}
which has real part
\begin{equation}
\frac{X}{X^2+(t-\gamma)^2} + \frac{\frac12 }{\frac14 + \gamma^2}.
\end{equation}

As in the polynomial case, we rescale:
\begin{align}
X=& x /\log t \cr
\Theta=&2\pi \theta /\log t .
\end{align}
Note that this is analogous to the rescaling in the
unitary case because $N\approx\log(t/2\pi)$.
We have
\begin{equation}\label{eqn:scaledx}
0= b^{**} +I -\frac12 \log t +\frac{2 x\log t}{x^2+\pi^2\theta^2}
  + x \log t \sum_{\gamma_*} \frac{1}{x^2+ 4\pi^2 \gamma_*^2 },
\end{equation}
where $b^{**}$ is a real constant, $I$ is purely imaginary, and
the sum is over the rescaled zeros $\gamma_* = \log t (t-\gamma)/2\pi$.
We follow the same procedure as in the $\UN$ case.
First multiply through~\eqref{eqn:scaledx} by $(x^2+\pi^2\theta^2)/\log t$ and
expand the final summand as a series in~$x$,
giving
\begin{equation}\label{eqn:xthetaalpha}
0=2x +(x^2+\pi^2 \theta^2)\left(-\frac12 + x\alpha_1-x^3\alpha_2+O(x^5)\right)+\text{ smaller terms},
\end{equation}
where
\begin{equation}
\alpha_j=\sum_{\gamma_*} \frac{1}{(4\pi^2 \gamma_*^2)^j }.
\end{equation}
Note that \eqref{eqn:xthetaalpha} has the same form as~\eqref{eqn:deltathetarelationN}.
Now write
$x=\beta_1 \pi^2 \theta^2 + \beta_2\pi^4\theta^4+O(\theta^6)$ and gather terms to
get
\begin{equation}
\left(2\beta_1-\frac12\right)\pi^2\theta^2
+
\left(2\beta_2-\frac12\beta_1^2+\beta_1\alpha\right) \pi^2\theta^4 + \text{ smaller terms} .
\end{equation}
Thus,
\begin{align}
\beta_1\sim\mathstrut& \frac14 \\
\beta_2 \sim\mathstrut& \frac{1}{64} -\frac18 \alpha_1,
\end{align}
which exactly corresponds to the $\UN$ case~in
Proposition~\ref{prop:delta}.  

The sum over zeros $\alpha_1$ is similar to the expression for $A_1$ in the
unitary case.  We can determine the expected value of the sum
if we assume $\mathrm{CUE}$ statistics for the zeros of the Riemann zeta-function,
restrict to having two closely spaced zeros, and find the one-level
density of the remaining zeros.  That calculation is in Section~\ref{sec:1level}.
In the notation of Lemma~\ref{lem:1level} we have
\begin{equation}
\langle \alpha_1 \rangle = \frac{1}{4\pi^2}
\int_{-\infty}^\infty \frac{1}{t^2} W_1^{(2,0)}(t)\, dt = \frac{1}{15}.
\end{equation}
Note that this is the same as the expected value of $A_1$.

Thus, assuming that the spacing of zeros of the Riemann zeta function has the
same distribution as the spacing of eigenvalues of random unitary matrices,
we find that the leading order behavior of zeros of $\zeta'$ near the
$\tfrac12$-line is the same as that of zeros of $\Lambda'$ near the unit circle.

%\subsection{Data for the zeta function}

\begin{figure}[htp]
\begin{center}
\scalebox{1.0}[1.0]{\includegraphics{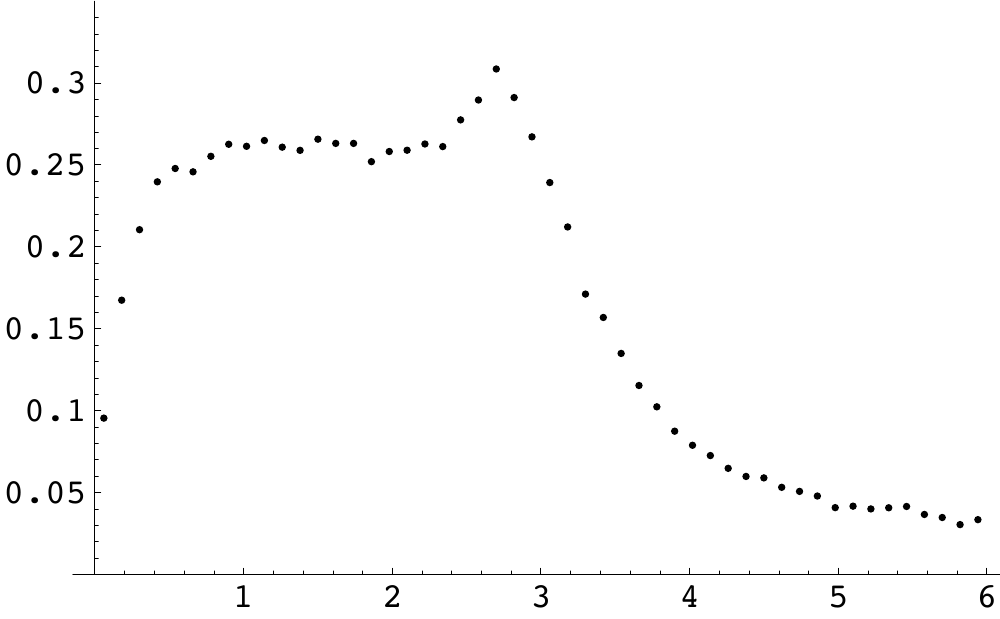}}
\caption{\sf
Normalized distribution of the real part of the zeros of $\zeta'(s)$.
Data is for the approximately 100000 zeros
with imaginary part in~$[10^6,10^6+60000]$.
} \label{fig:rezetaprime}
\end{center}
\end{figure}

In Figure~\ref{fig:rezetaprime} we show the empirical distribution of
the zeros of $\zeta'$ for $10^6 < t < 10^6 + 60000$.  The general shape
of the distribution shows a striking similarity with the zeros of
the derivative of characteristic polynomials of unitary matrices.

%Figure~\ref{fig:zetaIp12} compares the empirical cumulative distribution
%function of the normalized zeros of~$\zeta'$ with the predicted leading
%term.
%[[Do we want to eliminate that plot? --DF]]
%
%\begin{figure}[htp]
%\begin{center}
%\scalebox{0.75}[0.75]{\includegraphics{bestzetaguess12}}
%\caption{\sf
%Tail of $Ip(x)$ for the zeta-function near $10^6$
%compared to our prediction.
%} \label{fig:zetaIp12}
%\end{center}
%\end{figure}

\section{Calculation of the 1-level density}\label{sec:1level}

We prove the following
\begin{lemma}\label{lem:1level}
Fix $a,b>-1/2$.  If $t_1,t_2,\ldots,t_M$ are independently
  distributed with respect to the probability measure
  \begin{equation}
  \label{eqn:measurek}
  d\mu = d\mu^{(a,b)} =
  \frac1{C_{M,a,b}}\prod_{1\le n<m\le M}
        \left|e^{i t_n} - e^{i t_m}\right|^2
\prod_{1\le n\le M}
        \left|e^{i t_n} -1\right|^{2a}\left|e^{it_n}+1\right|^{2b}
dt_1 \cdots dt_{M} ,
\end{equation}
then the large-$M$ limiting (rescaled) 1-level density of the
normalized values $\tilde t_j = t_j M/2\pi $ is given by
\begin{equation}
  \label{eq:W-1}
  W_1^{(a,b)}(t)
  =t \frac{\pi^2}{2} \left(J_{a-\frac12}(\pi t)^2+J_{a+\frac12}(\pi t)^2\right)
  -a \pi J_{a-\frac12}(\pi t) J_{a+\frac12}(\pi t).
\end{equation}
\end{lemma}

The measure $d\mu^{(a,b)}$ is Haar measure on $U(M+a+b)$ restricted to those
matrices which have $a$ eigenvalues equal to~$1$ and $b$ eigenvalues equal to~$-1$.  Note that
$W_1^{(a,b)}$ is independent of~$b$.

\begin{proof}
  For fixed $a,b>-\frac12$ define $\omega(z)=|z-1|^{2a}|z+1|^{2b}$.
  Let $\{\phi_n\}_{n=0}^\infty$ be the sequence in $\C[x]$ uniquely
  determined by the following requirements:
  \begin{enumerate}
  \item For all $n\geq0$, $\phi_n$ is of degree~$n$ and has positive
    leading coefficient.
  \item For all $m,n\geq0$,
    \begin{equation}
      \label{eq:IPC}
      \langle\phi_m,\phi_n\rangle := \frac1{2\pi}
      \int_{-\pi}^{\pi}\phi_m(e^{it})\overline{\phi_n(e^{it})}
      w(e^{it})dt = \delta_{mn},
    \end{equation}
    where $\delta_{mn}$ is the Kronecker delta.
  \end{enumerate}
  Then $\{\phi_n\}$ is the sequence of normalized orthogonal polynomials on the
  unit circle with respect to the measure $d\nu(z)=(2\pi)^{-1}\omega(z)d\ell(z)$
  (where $d\ell$ is the arc-length element.)

 Let
 \begin{equation}
   \label{eq:K_M}
   K_M(z,w) := \sum_{n=0}^{M-1}\overline{\phi_n(z)}\phi_n(w)
 \end{equation}
 be the projection kernel onto polynomials of degree less than $M$
 with respect to the inner product~\eqref{eq:IPC}.  By the
 Gaudin-Mehta method, the probability measure~\eqref{eqn:measurek},
 when regarded as a measure on the unit circle, can be rewritten~as
 \begin{equation}
   \label{eq:dmu_k-GM}
   d\mu_k =\frac1{M!}
   \det_{1\leq j,k\leq M}(K_M(z_j,z_k))\prod_{j=1}^Md\nu(z_j)
 \end{equation}
 (note that the normalization constant $C_{M,k}$ is no longer needed).
 Then the $1$-level measure is $W_1^{(M)}(z)d\nu(z)$, where
 \begin{equation}
   \label{eq:W1_M}
   W_1^{(M)}(z) = K_M(z,z).
 \end{equation}
 (The normalization above is such that the total mass of the $1$-level
 measure is equal to $M$.)

 Let the ``dual'' $\phi^*$ of a polynomial
 $\phi(z)=c_nz^n+c_{n-1}z^{n-1}+\dots+c_1z+c_0$ of degree~$n$ be the
 polynomial
 \begin{equation}
   \label{eq:dual}
   \phi^*(z) = z^n\bar\phi(z^{-1})
   = \bar c_n + \bar c_{n-1}+\dots+\bar c_1z^{n-1}+\bar c_0z^n.
 \end{equation}
 Then we have the following formula of Szeg\H{o} for the projection
 kernel
(\cite{Sze39}, Theorem~11.4.2):
 \begin{equation}
   \label{eq:Szego-kernel}
   K_M(z,w) =
   \frac{\overline{\phi^*_M(z)}\phi^*_M(w) -
     \overline{\phi_M(z)}\phi_M(w)}
   {1-\bar z w}.
 \end{equation}
(This formula is analogous to the classical one of
   Christoffel and Darboux for the projection kernel of orthogonal
   polynomials on the line.)

 In view of~\eqref{eq:W1_M} and~\eqref{eq:Szego-kernel}, in order to
 find the rescaled limit of the $1$-level measure as $M\to\infty$ it
 will suffice to derive the asymptotic behavior of the orthogonal
 polynomials $\phi_n$ as $n\to\infty$ near the point $z=+1$.
 Theorem~\ref{thm:phi-vs-Jacobi} and formula~\eqref{eq:Hilb} below are
 the key ingredients, but first we need to introduce some notation.

 Denote by $\Pabn$ the classical Jacobi polynomials: they are
 orthogonal in the interval $[-1,1]$ with respect to the measure
 \begin{equation}
   \label{eq:w}
   w^{(a,b)}(x):=(1-x)^a(1+x)^b
 \end{equation}
 and are normalized as follows (\cite{Sze39}, Equation~4.3.3):
 \begin{equation}
   \label{eq:h_n}
   h_n^{(a,b)} := \int_{-1}^1|\Pabn(x)|^2w^{(a,b)}(x)\,dx =
   \frac{2^{a+b+1}}{2n+a+b+1}\,
   \frac{\Gamma(n+a+1)\Gamma(n+b+1)}{\Gamma(n+1)\Gamma(n+a+b+1)}.
 \end{equation}
 Let also
 \begin{align}
   \label{eq:hpm}
   h^+_n &= 2^{a+b}h_n^{(a+1/2,b+1/2)}, &
   h^-_n &= 2^{a+b}h_n^{(a-1/2,b-1/2)}.
 \end{align}

 \begin{theorem}\cite{Sze39}
   \label{thm:phi-vs-Jacobi}
   \begin{align}
     \notag
     z^{-n}\phi_{2n}(z) &=
     AP_n^{(a-1/2,b-1/2)}\left(\frac{z+z^{-1}}2\right)
     +B(z-z^{-1})P_{n-1}^{(a+1/2,b+1/2)}\left(\frac{z+z^{-1}}2\right) \\
     \label{eq:phi_e/o}
     z^{-n+1}\phi_{2n-1}(z) &=
     CP_n^{(a-1/2,b-1/2)}\left(\frac{z+z^{-1}}2\right)
     +D(z-z^{-1})P_{n-1}^{(a+1/2,b+1/2)}\left(\frac{z+z^{-1}}2\right).
   \end{align}
   Letting $c_n=(a+b)/(n+a+b)$, we have
   \begin{align*}
     A &= \sqrt{\frac\pi2\,\frac{1+c_n}{h_n^-}} &
     B &= \frac12\sqrt{\frac\pi2\,\frac{1-c_n}{h_{n-1}^+}} \\
     C &= \sqrt{\frac\pi2\,\frac{1-c_n}{h_n^-}} &
     D &= \frac12\sqrt{\frac\pi2\,\frac{1+c_n}{h_{n-1}^+}}.
   \end{align*}
 \end{theorem}
 Equation~\eqref{eq:phi_e/o} appears as~(11.5.4) in Szeg\H{o}'s book
 (except for an obvious typographical mistake therein.)  The constants
 $A,B,C,D$ can be easily found using Szeg\H{o}'s equation~(11.5.2)
 together with the fact that $\phi_{2n-1}$ is a polynomial, which
 forces the coefficient of $z^{-n}$ on the right-hand side of
 equation~\eqref{eq:phi_e/o} to vanish.  We omit the details.

 A meaningful rescaling of the $1$-level measure is achieved via the
 change of variables
 \begin{align}
   \label{eq:xi}
   t&=\frac{2\pi\xi}M.
 \end{align}
 Indeed, one finds that the limiting $1$-level measure (as
 $M\to\infty$) is $W_1(\xi)d\xi$, where
 \begin{align}
   \label{eq:K_inf}
   W_1(\xi)&=K_\infty(\xi,\xi), \\
   K_\infty(\xi,\eta) &= \lim_{M\to\infty}\frac1M
   K_M(e^{2\pi\xi/M},e^{2\pi\eta/M})\sqrt{\omega(e^{2\pi\xi/M})\omega(e^{2\pi\eta/M})}.
 \end{align}
(Actually, all limiting local correlations can be expressed
   in terms of the limiting kernel $K_\infty$, not just the $1$-level
   density.)

   It remains to compute $K_\infty(\xi,\eta)$.  It suffices to use
   equation~\eqref{eq:phi_e/o} in formula~\eqref{eq:Szego-kernel} and
   an asymptotic formula by Szeg\H{o}'s (see~\cite{Sze39},
   equation~(8.21.17)):
 \begin{equation}
   \label{eq:Hilb}
   \left(\sin\frac t2\right)^a\left(\cos\frac t2\right)^b\Pabn(\cos t)
   = N^{-a}\frac{\Gamma(n+a+1)}{n!}\sqrt{\frac t{\sin t}}J_a(Nt)
   + t^{a+2}O(n^a),
 \end{equation}
 valid for fixed $a>-1$, $b\in\R$, in the range $0<t\leq c/n$ for any
 fixed constant $c>0$, where $N=n+(a+b+1)/2$ and $J_a$ is a Bessel
 function of the first kind.

A tedious but straightforward computation finally gives:
\begin{align}
  \label{eq:K_infty}
  K_\infty(\xi,\eta) &=
  \frac\pi2 e^{i\pi(\eta-\xi)}\frac{\sqrt{\xi\eta}}{\xi-\eta}
  \left(J_{a+1/2}(\pi\xi)J_{a-1/2}(\pi\eta)
    - J_{a-1/2}(\pi\xi)J_{a+1/2}(\pi\eta)\right)\\
  K_\infty(\xi,\xi) &=
  \frac\pi2\left\{ \pi\xi[J_{a+1/2}(\pi\xi)^2+J_{a-1/2}(\pi\xi)^2]
    - 2a J_{a+1/2}(\pi\xi)J_{a-1/2}(\pi\xi) \right\}.
\end{align}
\begin{remark}
  The kernel $K_{\infty}(\xi,\eta)$ can be used to compute $n$-level
  correlations and spacing statistics through (matrix or operator)
  determinants (see~\cite{TW98} for an explanation and proof of these
  applications). Incidentally, for the purposes of evaluating such
  determinants, the factor $e^{i\pi(\eta-\xi)}$ may be suppressed
  in~\eqref{eq:K_infty} (this is tantamount to conjugating the
  corresponding integral operator by a unitary transformation).
\end{remark}

\end{proof}

\section{Appendix: The Domain of Analyticity of the Root of the Derivative}
\label{sec:dom-analyt-root}
\begin{center}
{by Eduardo Due\~nez}
\end{center}

Define as before $e(x):=e^{2\pi ix}$ for $x$ real. For any fixed
$N\geq3$ define $e_N(x):=e(x/N)=e^{2\pi ix/N}$.  Consider $N$ unit
complex numbers
\begin{equation}
  e_N(\theta_0),e_N(-\theta_0),
  e_N(\theta_1),e_N(\theta_2),\dots,e_N(\theta_{N-2}),
  \label{eq:roots-f}
\end{equation}
where
\begin{equation}\label{eq:pyram}
  0\leq\theta_0\leq\theta_j\leq N-\theta_0,\qquad j=1,2,\dots,N-2.
\end{equation}
(I.~e., the open arc centered at $1$ of the unit circle $|z|=1$ having
endpoints $e_N(\pm\theta_0)$ contains none of the remaining $N-2$
numbers). The inequalities~(\ref{eq:pyram}) define a ``pyramid''
$\mathcal{P}_N$ contained in the cube $[0,N]^{N-1}$. For fixed
$0<T\leq N/2$, denote by $\mathcal{P}_N^{(T)}$ the closed truncation
of $\mathcal{P}_N$ at height~$T$. Then $\mathcal{P}_N^{(T)}$ is defined by
the inequalities
\begin{align}
  \theta_0&\leq\theta_j\leq N-\theta_0,\qquad j=1,2,\dots,N-2;\notag\\
  \label{eq:pyram-T}0&\leq\theta_0\leq T.
\end{align}
For notational convenience we will set
$\Theta=(\theta_1,\theta_2,\dots,\theta_{N-2})$ and
$\Theta_0=(\theta_0;\Theta)$. Finally, let
\begin{equation}\label{eq:9}
  f(z)=f(\Theta_0;z)=
  (z-e_N(\theta_0))(z-e_N(-\theta_0))\prod_{j=1}^{N-2}(z-e_N(\theta_j))
\end{equation}
be the monic polynomial of degree~$N$ with roots~(\ref{eq:roots-f}).

In this section we prove the following result.
\begin{proposition}\label{prop:doma-analyt-root}
  With the above notation, for every
  $T<\frac1\pi$ there exists $N(T)$ such that for all $N\geq
  N(T)$ and for each $\Theta_0$ in $\mathcal{P}_N^{(T)}$, the
  derivative $f'(z)=\frac{\partial}{\partial z}f(\Theta_0;z)$ of $f$
  has a unique root $z'=\varsigma(\Theta_0)$ in the open disk with
  diameter $[e_N(-\theta_0),e_N(\theta_0)]$. Moreover, $\varsigma$ is an
  analytic function of $\Theta_0$ in the interior of
  $\mathcal{P}_N^{(T)}$.
\end{proposition}
We remark that, by the Gauss-Lucas theorem, the root $z'$ alluded to
in Proposition~\ref{prop:doma-analyt-root} must lie on or to the
left of the vertical diameter~$[e_N(-\theta_0),e_N(\theta_0)]$.
\begin{proof}
  Fix $T<1/\pi$ and
  consider $N$ as a parameter ($N\geq3$) for the time being. Let
  \begin{equation}
    \label{eq:4}
    F(z)=F(\Theta;z):=\prod_{j=1}^{N-2}(z-e_N(\theta_j)),
  \end{equation}
  so
  \begin{equation}
    f(z)=(z-e_N(\theta_0))(z-e_N(-\theta_0))F(z)\label{eq:6}
  \end{equation}
  and
  \begin{equation}
    \label{eq:7}
    \frac{f'}{f}(z) = g(z) + L(\Theta;z),
  \end{equation}
  where
  \begin{align}
    g(z)&:=\frac1{z-e_N(\theta_0)}+\frac1{z-e_N(-\theta_0)},\label{eq:12}\\
    L(\Theta;z) &:= \frac{F'}{F}(z)
    = \sum_{j=1}^{N-2}\frac1{z-e_N(\theta_j)}. \label{eq:8}
  \end{align}

  Define $c_N(x):=\cos(2\pi x/N)$ and $s_N(x):=\sin(2\pi x/N)$, so
  $e_N(x)=c_N(x)+is_N(x)$. Parametrize the boundary of the disk with
  diameter $e_N(\pm\theta_0)$ (i.~e., the disk
  $|z-c_N(\theta_0)|\leq s_N(\theta_0)$) as:
  \begin{equation}
    z(\phi)=c_N(\theta_0)+ie^{i\phi}s_N(\theta_0).
\label{eq:10}
  \end{equation}
  Since roots of $f$ occur at $e_N(\pm\theta_0)$ whenever $f$ has
  multiple roots there, it is best to work instead with a
  slightly deformed contour $\mathcal{C}$ obtained as the boundary of
  the ``twice bitten'' disk
  \begin{align*}
    |z-c_N(\theta_0)|&\leq s_N(\theta_0)\\
    |z-e_N(\theta_0)|&\geq \epsilon \\
    |z-e_N(-\theta_0)|&\geq \epsilon.
  \end{align*}
  missing tiny $\epsilon$-neighborhoods of the points
  $e_N(\pm\theta_0)$. We still assume that $\mathcal{C}$ is parametrized
  by~(\ref{eq:10}), except for $\phi$ in a small $\delta$-neighborhood
  of any multiple of~$\pi$. (We will not write down the exact
  parametrization of $\mathcal{C}$ for such values of $\phi$ since
  its precise form  will  not be needed.)

  Write $\mathcal{C}$ as a union of four pieces: $\mathcal{C}
  =\mathcal{C}_\text{d}\cup\mathcal{C}_\text{r}\cup\mathcal{C}_\text{u}
  \cup\mathcal{C}_\ell$ (down, right, up, left), where
  \begin{align*}
    \mathcal{C}_\text{d}&=\{z(\phi):\pi-\delta\le\phi\le\pi+\delta\},\\
    \mathcal{C}_\text{r}&=\{z(\phi):-\pi+\delta\le\phi\le-\delta\},\\
    \mathcal{C}_\text{u}&=\{z(\phi):-\delta\le\phi\le\delta\},\\
    \mathcal{C}_\ell&=\{z(\phi):\delta\le\phi\le\pi-\delta\}.
  \end{align*}

  The assumed inequalities~(\ref{eq:pyram-T}) ensure that
  there are no zeros of $f$ anywhere on $\mathcal{C}$ (and,
  \textit{a fortiori,} no multiple zeros); hence,
  every
  zero of $f'$ inside $\mathcal{C}$ is a zero of $f'/f$ with the
  same multiplicity. It suffices to ensure that
  that (for all sufficiently small $\epsilon$)
  the contour $\mathcal{C}$ encloses exactly one zero of~$f'/f$. Note
  that, since $f(\Theta_0;z)$ has no zeros on or inside $\mathcal{C}$,
  $f'/f$ has no poles there either.  By the Argument Principle, the
  claim in Proposition~\ref{prop:doma-analyt-root} regarding the
  uniqueness
  of the zero $z'$ of $f'$ is equivalent to showing that the image
  $\mathcal{D}$ of
  $\mathcal{C}$ under
  $f'/f$ has index $1$ about the origin for large $N$ and all
  sufficiently small~$\epsilon$.

  We let
  $w(\phi):=\frac{f'}{f}(z(\phi))$ and
  denote the images of the four pieces $\mathcal{C}_\text{d},
  \mathcal{C}_\text{r},\mathcal{C}_\text{u},\mathcal{C}_\ell$ of
  $\mathcal{C}$ under $f'/f$ by
  $\mathcal{D}_\text{d},\mathcal{D}_\text{r},\mathcal{D}_\text{u},
  \mathcal{D}_\ell$.

  The idea of the proof is very simple. As we shall show, the dominant
  part of the logarithmic derivative $f'/f(z)$ (for $z$ in
  $\mathcal{C}$) is $g(z)$.
  The image of $\mathcal{C}$ under $g(z)$ is easy to describe
  explicitly; indeed, it is a curve $\mathcal{E}$ having index $1$
  about the
  origin. The curve $\mathcal{D}$ can be regarded as a perturbation of
  $\mathcal{E}$. We show that, if $H$ and $\epsilon$ are sufficiently
  small, then the remaining terms $1/(z-e_N(\theta_j))$,
  $j=1,2,\dots,N-2$, are small enough to keep the index of
  $\mathcal{D}$ equal to that of~$\mathcal{E}$.

  \begin{figure}[h]
    \centering
    \scalebox{.4}[.4]{\includegraphics{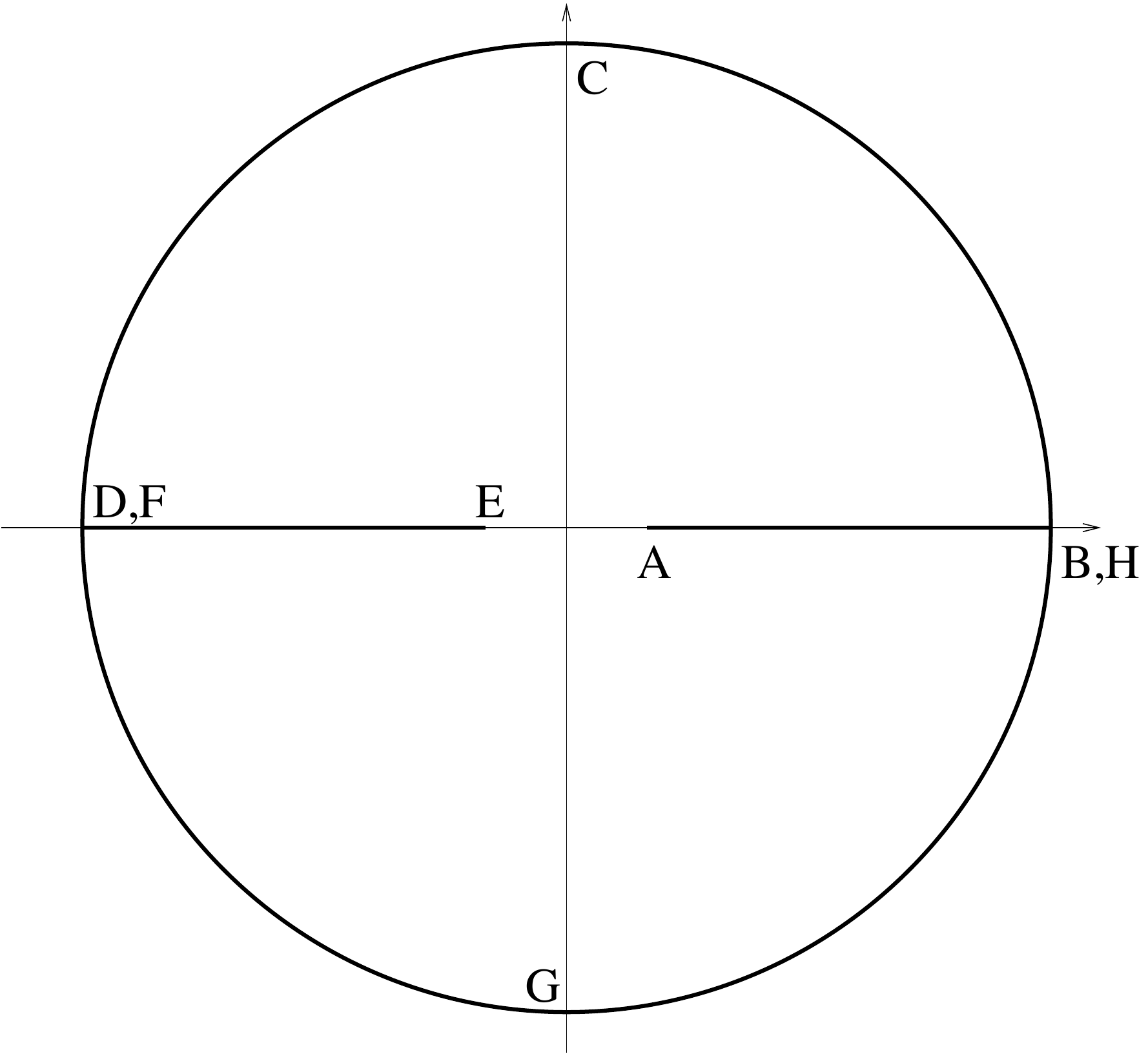}}
    \caption{The curve $\mathcal{E}$.}
    \label{fig:contourE}
  \end{figure}
  We denote by
  $\mathcal{E}_{\text{d}},\mathcal{E}_{\text{r}},\mathcal{E}_{\text{u}},
  \mathcal{E}_\ell$ the images under $g$ of the respective parts of
  $\mathcal{C}$ (figure~\ref{fig:contourE}).

  \begin{itemize}
  \item $\mathcal{E}_{\text{r}}$ is parametrized as
    \begin{equation}
      g(z(\phi))=\frac1{2s_N(\theta_0)\sin\phi},\qquad
      -\pi+\delta\leq\phi\leq-\delta.\label{eq:11}
    \end{equation}
    Thus, $\mathcal{E}_{\text{r}}$ is a twice-traversed straight
    segment starting at the faraway point
    $H=g(z(-\pi+\delta))=1/(2s(\theta_0)\sin\delta)$, moving
    leftward to the point $A=g(z(0))=1/(2s(\theta_0))$ and retracing
    itself back to~$B=g(z(-\delta))$ (here $B=H$).\vspace{1ex}
  \item $\mathcal{E}_{\text{u}}$ is a large arc $\overarc{BCD}$ on the
    upper half-plane. ($\mathcal{E}_{\text{u}}$ is essentially a
    semicircle, for $\epsilon$ small. This is easily seen from the
    fact that the dominant term in $f'/f(z)$ for $z$ very close to
    $e_N(\theta_0)$ is $1/(z-e_N(\theta_0))$, namely an inversion with
    center
    $e_N(\theta_0)$, taking the tiny (almost) semicircle
    $\mathcal{C}_{\text{u}}$ to a huge (almost)
    semicircle~$\mathcal{E}_{\text{u}}$.)  $\mathcal{E}_{\text{u}}$
    escapes any bounded region of the plane as $\epsilon$ approaches
    zero.
  \item $\mathcal{E}_{\ell}$ and $\mathcal{E}_{\text{d}}$ are obtained
    from $\mathcal{E}_{\text{r}}$ and $\mathcal{E}_{\text{u}}$ through
    central symmetry with respect to the origin.
  \end{itemize}
  It is clear that $\mathcal{E}$ has index $1$ about the origin for
  all sufficiently small $\epsilon$.

  For notational convenience, denote by $\tilde A=w(-\pi/2),\tilde
  B=w(-\delta),\dots,\tilde H=w(-\pi+\delta)$ the points on
  $\mathcal{D}$ analogous to $A,B,\dots,H$ on~$\mathcal{E}$.

  We claim that if $N$ is large enough and
  $\Theta_0\in\mathcal{P}_N^{(T)}$, then
  the index of $\mathcal{D}$ about the origin is also~$1$ for all
  sufficiently small~$\epsilon$ and all $N\geq N(T)$. Note that
  $\epsilon$ is allowed to depend on $\Theta_0$, whereas the lower
  bound $N(T)$ for $N$ depends only on~$T$.

  For notational simplicity, we set $\theta_{N-1}:=-\theta_0$. Write
  \begin{equation*}
    \frac{f'}{f}(z) = \frac{m}{z-e_N(\theta_0)}
    +
    \underbrace{\sum_{j=m}^{N-1}\frac1{z-e_N(\theta_j)}}_{\displaystyle
      R(\Theta_0;z)},
  \end{equation*}
  where $m\geq1$ is the multiplicity of the zero $e_N(\theta_0)$ of
  $f(z)$ and $\theta_j$, $j=m,\dots,N-1$, are the remaining zeros
  of~$f(z)$.

  Let $\epsilon_0>0$ be the minimum of the distances from
  $e_N(\theta_0)$
  to the other $e_N(\theta_j)$, $m\leq j\leq N-1$.
  Now let
  $\epsilon=\epsilon_0/N$. If $|z-e_N(\theta_0)|=\epsilon$, then
  $|z-e_N(\theta_j)|\geq\epsilon_0$ ($m\leq j\leq N-1$), so
  $|R(\Theta_0;z)|\leq(N-m)/\epsilon_0$. Moreover,
  $m/(z-e_N(\theta_0))$ maps $\mathcal{C}_{\text{u}}$ into (part of) the
  upper half-circle $\mathcal{S}$: $|w|=m/\epsilon$, $\Re
  w>0$. Therefore, $\mathcal{D}_{\text{u}}$ will be supported on the
  $\frac{N-1}{\epsilon_0}$-neighborhood $\mathcal{T}$
  of~$\mathcal{S}$. Since $\mathcal{T}$ intersects the imaginary axis
  $\Re w=0$ on the interval
  $i\left[\frac m\epsilon-\frac{N-m}{\epsilon_0},
  \frac m\epsilon+\frac{N-m}{\epsilon_0}\right]
  =i\left[\frac{m+N(m-1)}{\epsilon_0},\frac{N+m(N-1)}{\epsilon_0}\right]$,
  our claim
  that $\mathcal{D}_{\text{u}}$ only intersects the positive imaginary
  axis is proved. The claim that $\mathcal{D}_\ell$ only intersects
  the negative imaginary axis for suitably small $\epsilon$ is proved
  along identical lines.

  In order to prove the claim that $\mathcal{D}$ has index~$1$ about
  the origin, it suffices now to show that $\mathcal{D}_{\text{r}}$ is
  contained in the right half-plane $\Re(w)>0$ and
  $\mathcal{D}_{\ell}$ in the left half-plane $\Re(w)<0$. The truth of
  said statement for $\mathcal{D}_{\text{r}}$ is obvious, since each
  of the terms $1/(z-e_N(\theta_j))$ ($j=0,\dots,N-1$) in $f'/f(z)$
  (cf., equations~(\ref{eq:7})--(\ref{eq:8})) has positive real part
  for $z$ on $\mathcal{C}_{\text{r}}$.

  Now note that
  \begin{equation}
    \label{eq:13}
    -g(z(\phi))=\frac1{2s_N(\theta_0)\sin\phi},\qquad
          \delta\leq\phi\leq\pi-\delta,
  \end{equation}
  is real and positive; we anticipate it to be the
  dominant term of the logarithmic derivative $-f'/f(z)$. We will show
  that the real parts of each of the terms
  $\frac1{z(\phi)-e_N(\theta_j)}$ (for
  $\delta\leq\phi\leq\pi-\delta$ and $j=1,\dots,N-2$) are bounded above
  by a suitably small multiple of $-g(z(\phi))$.

  \begin{lemma}\label{lem:bound-reciprocals}
    For all $z$ with $|z|<1$ we have
    \begin{equation}
      \label{eq:17}
      \max_{|\zeta|=1}\Re\frac{1}{z-\zeta} = \frac{1-\Re(z)}{1-|z|^2}.
    \end{equation}
  \end{lemma}
  This lemma can be proved by rewriting the equation $|\zeta|=1$
  in the form
  \begin{equation*}
    \left|w+\frac{\bar z}{1-|z|^2}\right| = \frac1{1-|z|^2},
  \end{equation*}
  in terms of the variable
  \begin{equation*}
    w = \frac{1}{z-\zeta},
  \end{equation*}
  whence the result follows trivially.

  \begin{lemma}\label{lem:bound-eta}
    For $N\geq3$, $0\leq\theta_0\leq1$, $0<\phi<\pi$ and all $\psi$ we
    have
    \begin{equation}
      \label{eq:18}
      \Re\frac{1}{z(\phi)-e_N(\psi)}
      \leq -\left(\frac{\pi\theta_0}N
        + 2\pi^2(7+4\sqrt3)\frac{\theta_0^2}{N^2}\right)g(z(\phi)),
    \end{equation}
  \end{lemma}
  First, since $-g(z(\phi))>0$ for $0<\phi<\pi$ (cf.,
  equation~(\ref{eq:13})), it suffices to prove
  the upper bound
  \begin{equation}\label{eq:20}
    \left(\frac{\pi\theta_0}N
      + 2\pi^2(7+4\sqrt3)\frac{\theta_0^2}{N^2}\right)
  \end{equation}
  for the quantity
  \begin{equation}
    \label{eq:19}
    h(\theta_0;\phi,\psi)
    := \frac{-1}{g(z(\phi))}\Re\frac{1}{z(\phi)-e_N(\psi)}
    = s_N(\theta_0)\sin\phi\Re\frac{1}{z(\phi)-e_N(\psi)}.
  \end{equation}
  By Lemma~\ref{lem:bound-reciprocals},
  \begin{equation*}
    h(\theta_0;\phi,\psi) \leq
    \left(\frac{1-\Re(z(\phi))}{1-|z(\phi)|^2}\right) s_N(\theta_0)\sin\phi.
  \end{equation*}
  From Equation~(\ref{eq:10}) it quickly follows that
  \begin{align*}
    1-|z(\phi)|^2 &= s_N(2\theta_0)\sin\phi,\quad\text{and}\\
    \Re(z(\phi)) &= c_N(\theta_0)-s_N(\theta_0)\sin\phi.
  \end{align*}
  Thus,
  \begin{align}\label{eq:21}
    h(\theta_0;\phi,\psi)
    &\leq \left(\frac{1-c_N(\theta_0)+s_N(\theta_0)\sin\phi}
      {s_N(2\theta_0)\sin\phi}\right)s_N(\theta_0)\sin\phi
    =\frac{1-c_N(\theta_0)+s_N(\theta_0)\sin\phi}{2c_N(\theta_0)}\notag\\
    &= \frac12\tan\left(\frac{2\pi\theta_0}N\right)\sin\phi
    +\frac12\left(\sec\left(\frac{2\pi\theta_0}N\right)-1\right)\notag\\
    &\leq \frac12\tan\left(\frac{2\pi\theta_0}N\right)
    +\frac12\left(\sec\left(\frac{2\pi\theta_0}N\right)-1\right).
  \end{align}
  Let us now use the Taylor formul\ae\ with remainder
  \begin{align*}
    |\tan\theta-\theta| &\leq \frac{\theta^2}2 \sup_{|\vartheta|\leq\theta}|\tan''\vartheta|\\
    |\sec\theta-1| &\leq \frac{\theta^2}2 \sup_{|\vartheta|\leq\theta}|\sec''\vartheta|
  \end{align*}
  with $\theta:=2\pi\theta_0/N$.
  Both $|\tan''\vartheta|$ and
  $|\sec''\vartheta|$ are even functions of $\vartheta$, increasing
  with $|\vartheta|$,
  so their respective suprema are bounded by $|\tan''(2\pi/3)|=14$ and
  $|\sec''(2\pi/3)|=8\sqrt3$ (since
  $\vartheta=\theta=2\pi\theta_0/N\leq2\pi/3$, by the assumptions
  $\theta_0\leq1$ and $N\geq 3$). The inequalities
  \begin{align*}
    \tan\theta &\leq \theta + 7\theta^2,\quad\text{and}\\
    \sec\theta-1 &\leq 4\sqrt3\,\theta^2,
  \end{align*}
  follow immediately. These inequalities together with~(\ref{eq:21})
  complete the proof of Lemma~\ref{lem:bound-eta} upon setting
  $\theta=2\pi\theta_0/N$.

  To complete the proof that $\mathcal{D}_\ell$ is contained in $\Re
  w<0$, it remains to show that
  \begin{equation}
    \label{eq:22}
    -\Re\left(\frac{f'}{f}(z(\phi))\right)>0
    \qquad\text{for $\delta\leq\phi\leq\pi-\delta$.}
  \end{equation}
  But
  \begin{align*}
    -\Re\left(\frac{f'}{f}(z(\phi))\right) &=
    -g(z(\phi)) -
    \sum_{j=1}^{N-2}\Re\frac1{z(\phi)-e_N(\theta_j)}
    \qquad\text{since $-g(z(\phi))>0$} \\
    &\geq -g(z(\phi)) \left[1-(N-2)\left(\frac{\pi\theta_0}N
        + 2\pi^2(7+4\sqrt3)\frac{\theta_0^2}{N^2}\right)\right]
    \qquad\text{by Lemma~\ref{lem:bound-eta}}\\
    &= -g(z(\phi))
    \left[1-\pi\theta_0 + U\frac{\theta_0}{N} + V\frac{\theta_0^2}{N} + W\frac{\theta_0^2}{N^2}\right],
  \end{align*}
  say, for suitable absolute constants $U,V,W$. As $N\to\infty$, the last
  bracketed quantity above has the limit $1-\pi\theta_0$. As long as
  $|\theta_0|\leq T<\frac1\pi$, there will exist $N(T)$ such that said quantity
  is positive for $N\geq N(T)$. This completes the proof of~(\ref{eq:22}) and of
  the uniqueness of the desired root $z'=\varsigma(\Theta_0)$ of~$f'$.  The
  analyticity of $\varsigma$ as a function of $\Theta_0$ follows from the joint
  analyticity of the function $f(\Theta_0;z)$ in the variables $\Theta_0$ and
  $z$ together with the formula
  \begin{equation}
    \label{eq:23}
    \varsigma(\Theta_0) =
    \frac1{2\pi i}\oint_{\mathcal{C}}
    \frac{z f''(\Theta_0;z)}{f'(\Theta_0;z)}dz.
  \end{equation}
  While $\mathcal{C}$ depends on the choice of a fixed $\epsilon$, it
  is clear that $\epsilon$ can be chosen so that
  $z'=\varsigma(\Theta_0)$ is uniquely
  defined by the integral~(\ref{eq:23}) for all $\Theta_0$ in any
  desired
  compact subset of the interior of $\mathcal{P}_N^{(T)}$. This is
  enough to ensure that $\varsigma$ is analytic in the whole interior
  of $\mathcal{P}_N^{(T)}$ and concludes the proof of Proposition~8.1.
 \end{proof}

\end{document}